\input ./style/arxiv-general.cfg
\documentclass[aop,MSNbibl,nameyear,dvips]{arximspdf}
\makeatletter
   \@ifpackageloaded{graphicx}{}{\usepackage{graphicx}}
\makeatother
\usepackage{accents}

\doi{10.1214/15-AOP1032}
\volume{44}
\issue{4}
\pubyear{2016}
\firstpage{2694}
\lastpage{2725}
\docsubty{FLA}

\makeatletter
\def\sfrac#1#2{#1/#2}

\def\afrac#1#2{#1/(#2)}
\def\vafrac#1#2{(#1)/(#2)}

\renewcommand{\mathring}[1]{\accentset{\circ}{#1}}
\newcommand{\rrvert}{\vert}
\newcommand{\llvert}{\vert}
\newcommand{\eqref}[1]{(\ref{#1})}
\newtheorem{teo}{Theorem}[section]
\newtheorem{prop}[teo]{Proposition}
\newtheorem{cor}[teo]{Corollary}
\newproclaim{defin}[teo]{Definition}
\newtheorem{lem}[teo]{Lemma}
\newproclaim{rem}{Remark}[section]

\newcommand{\Diam}{\operatorname{Diam}}
\makeatother

\begin{document}
\begin{frontmatter}

\title{Improper Poisson line process as SIRSN in~any~dimension}
\runtitle{Random lines yield a SIRSN}

\begin{aug}
\author[A]{\fnms{Jonas}~\snm{Kahn}\corref{}\ead[label=e1]{jonas.kahn@math.univ-lille1.fr}}
\runauthor{J. Kahn}
\affiliation{Universit\'{e} de Lille 1, CNRS}
\address[A]{Laboratoire Paul Painlev\'{e} (UMR 8524)\\
Universit\'{e} de Lille 1, CNRS\\
Cit\'{e} Scientifique---B\^{a}t. M2\\
59655 Villeneuve d'Ascq Cedex\\
France\\
\printead{e1}}
\end{aug}

%
\received{\smonth{3} \syear{2015}}
%
\revised{\smonth{4} \syear{2015}}

%
\begin{abstract}
Aldous has introduced a notion of scale-invariant random spatial
network (SIRSN)
as a mathematical formalization of road networks. Intuitively, those
are random
processes that assign a route between each pair of points in Euclidean space,
while being invariant under rotation, translation, and change of scale,
and such that the routes are not too long and mainly lie on ``main roads''.

The only known example was somewhat artificial since invariance had to be
added using randomization at the end of the construction. We prove that the
network of geodesics in the random metric space generated by a Poisson line
process marked by speeds according to a power law is a SIRSN, in any dimension.

Along the way, we establish bounds comparing Euclidean balls and balls for
the random metric space. We also prove that in dimension more than two, the
geodesics have ``many directions'' near each point where they are not straight.
\end{abstract}

%
\end{frontmatter}

\section{Introduction}
\label{intro}

Scale-invariant random spatial networks (SIRSNs) are a class of random
networks defined as a route between each pair of points, with three
types of properties. First, invariance properties, second, guarantees
on mean lengths of routes---in the Euclidean metric---and third
guarantees on the mean length of intersection of a suitably truncated
version of the network with a given compact set. It turns out that
these conditions are enough to deeply constrain the network. For
example, all SIRSNs have singly-infinite paths for which any subset is
included in a route, but no such paths are doubly-infinite.

The only known example is the binary hierarchy model, in two
dimensions. It consists of minimum-time paths on a dyadic grid where
speed depends on the two-valuation. Invariance is obtained by a
randomization construction. The latter feature is somewhat
unsatisfying: the model itself is invariant, but realisations exhibit
long-range dependence: observation of a small region gives much
information on the network everywhere.

A more ``natural'' candidate for a SIRSN is therefore the \emph{Poisson
line process model}. Intuitively, lines are thrown uniformly at random
in $\mathbb{R} ^d$, and marked with random speed limits. Slower lines
are dense in $\mathbb{R} ^d$. Then the route between two points is the
minimum-time path made of segments of these lines. Remarkably, even in
dimension $d\geq3$, when random lines almost surely do not intersect,
such paths exist. The whole construction is invariant. The aim of this
paper is to show that the Poisson line process is indeed a SIRSN for
all $d\geq2$.

Historically, \citet{Aldous2014} introduced the notion of a SIRSN,
and proved a number of their properties, including those mentioned in
the first paragraphs of this \hyperref[intro]{Introduction}. \citet{Aldous.Ganesan}
give a verbal description. The motivation was twofold.

First, \citet{Aldous.Kendall} had proved that it was possible to
build a road network connecting a prescribed set of points that both
had routes almost as short as the segments between each pair of points,
that is, the corresponding Euclidean geodesics, and total road length
almost as short as the Steiner tree, that is, the shortest possible
connecting network. However, the network was less efficient at small
scales. Thanks to their invariance properties, SIRSNs have the same
efficiency at all scales. It turns out that there is a trade-off: for a
SIRSN, there is a lower bound on the expected total length of the
network, which is decreasing in the expected route length between two points.

The second motivation was to give a mathematical abstraction of road
networks and maps, in particular online maps as they are used today.
Namely, we may change viewpoint, zoom in or out, and the appearance
changes little, as smaller roads are shown and hidden. SIRSNs are then
defined as statistically invariant under translation, rotation and
change of scale. Moreover, we are less interested in the roads than in
the routes: how do we drive from $A$ to $B$? SIRSNs are thus defined by
giving routes only, namely unique routes connecting pairs of points. It
turns out that a notion of ``main roads'' at any scale can be built
from the network of routes itself. To wit, the network of main roads at
scale $r$ would be the network of routes deprived of balls of radius
$r$ around their endpoints. It is finite in every compact. Similarly,
\citet{Kalapala.} have shown that a number of statistics of real
road networks do not depend on scale.

\citet{Aldous2014}proved that the binary hierarchy model was a
SIRSN, and suggested two other possible models for SIRSN, one of which
is the Poisson line process model. \citet{Kendall2014} has then
proved important properties of the Poisson line process with
appropriate speeds: it does yield a random metric space and this space
is a geodesic space. Moreover, in dimension two, the geodesics are
almost everywhere unique, the geodesics are locally of finite
mean-length, and the subnetwork obtained from the routes connecting
points of an independent Poisson point process has finite length in a
compact set. The latter properties establish a ``pre-SIRSN'' result,
but fall short of the full definition.

As a candidate for a SIRSN, Poisson line process model then fall in a
large category: networks derived from geodesic spaces. Indeed, from any
geodesic space, we may build a spatial network by associating to any
two points the geodesic(s) connecting them. It is not obvious how one
might determine when such a network is a SIRSN.

In Section~\ref{notations}, we give a precise definition of a SIRSN and
of Poisson line processes. We also present other notation and
definitions, and recall some known results, in particular that the
Poisson line process with speed limits yields a random metric on
$\mathbb{R} ^d$. We then give a few basic properties of $\Pi$-paths,
that is paths in this metric space. In Section~\ref{time_diam}, we give
sharp bounds on the random diameter for this metric of a Euclidean
ball, with a few generalizations. These estimates will be an important
tool in several subsequent proofs. In Section~\ref{unique}, we prove
that geodesics between a given pair of points are almost surely unique,
in any dimension. Lemma~\ref{tour} will play a central role: we
introduce the notion of ``many directions'', and the lemma states that
geodesics have many directions at relevant points. A consequence is
that any path using the same lines as a geodesic will contain these
points. We will then conclude by noticing that geodesics between the
same pair of points almost surely use the same lines (Lemma~\ref
{sameseg}). In Section~\ref{expect_geo}, we prove that geodesics have
finite mean Euclidean length. Alternatively, we may see the result as
supplying a stochastic control of the Euclidean diameter of balls for
the metric generated by the Poisson line process. Section~\ref{ldn}
contains the last and most important component of the proof that the
Poisson line process generates a SIRSN. Intuitively, we establish a
sharp control of the total length of the intersection of all infinitely
many geodesics minus a ball around each endpoint, with a given ball.
This corresponds to the fact that all these geodesics coalesce before
hitting the ball and split after leaving it. Bounds are given using the
pigeon-hole principle and the fact that relevant geodesics must use the
few fast lines that hit the ball. Finally, Section~\ref{conclusion}
suggests a few potential directions of future inquiry.

\section{Notation, definitions, basic properties}
\label{notations}

We follow \citeauthor{Kendall2014}'s (\citeyear{Kendall2014}) notation
whenever possible.


We write $B(x, r)$ for the ball with center $x$ and radius $r$.

We first give the precise definition of a SIRSN. Suppose that $\Lambda
$ is an atom-free measure on a measurable space $(\mathcal{X}, \mathcal
{B} ) $. Recall that a Poisson point process of intensity measure
$\Lambda$ is a random set of points such that for any $B \in\mathcal
{B} $, the number of points $N(B)$ in $B$ is a Poisson variable with
intensity $\Lambda(B)$, and the number of points $N(B_i)$ in disjoint
$B_i$ are independent. Then a SIRSN is defined as follows.
\begin{defin}
\label{defSIRSN}
A SIRSN is a process that associates to any two points $x_1$ and $x_2$
in $\mathbb{R} ^d$ random routes such that:
\begin{longlist}[2.]
\item[1.]  Between two specified points $x_1$ and $x_2$,
there is almost surely a unique route $\mathcal{R} (x_1, x_2) = \mathcal
{R} (x_2, x_1)$. It is a finite-length path connecting $x_1$ and $x_2$.

\item[2.] For a finite number of points $x_1, \dots, x_k$
in $\mathbb{R} ^d$, consider the subnetwork $\mathcal{N} (x_1, \dots,
x_k)$ formed by the random routes connecting all $x_i$ and $x_j$. Then
$\mathcal{N} (x_1, \dots, x_k)$ is statistically invariant under
translation, rotation and change of scale: if $\mathfrak{R}$ is a
Euclidean similarity of $\mathbb{R} ^d$, then $\mathcal{N} (\mathfrak
{R}(x_1), \dots, \mathfrak{R}(x_k) )$ has the same distribution as
$\mathfrak{R}\mathcal{N} (x_1, \dots, x_k)$.

\item[3.] Let $D_1$ be the length of a route between two
points at unit distance. Then $\mathbb{E} [D_1 ] < \infty$.

\item[4.]  Let $ \{ \Xi_n, n\in\mathbb{N} ^* \} $ be a
collection of Poisson processes on $\mathbb{R} ^d$ with intensity $n$
times Lebesgue, all independent from the SIRSN. Suppose they are
coupled so that $\Xi_n \subset\Xi_{n+1}$. Write $\Xi= \bigcup_{n\in
\mathbb{N} ^*} \Xi_n$. Then the \emph{intensity} (mean length per unit
area) $p(1)$ of the following long-distance network is finite:
\[
\bigcup_{x_1, x_2 \in\Xi} \bigl( \mathcal{R} (x_1,
x_2) \bigr) \setminus\bigl(B(x_1, 1) \cup
B(x_2, 1)\bigr).
\]
\end{longlist}
\end{defin}
Note that \citet{Kendall2014} uses more conditions in his
definition, but the missing properties are implied by property~4. They were useful to define weaker variants of a SIRSN.

The use of Poisson processes in property~4 makes it look
slightly complex, but this is a technical shortcut: it allows us to
study the network through only countably many routes. Morally, we would
like property~4 to hold true if the long-distance network
was defined as the union of all routes between all pairs of points of
$\mathbb{R} ^d$, minus the balls around the endpoints. But there would
be uncountably many routes, and it would be harder to work with.

We now turn to the definition of the improper Poisson line process. We
first need a measure on lines. More details on this kind of process may
be found in the book by \citet{stocGeom}.

Let $\mathcal{L} ^d$ be the space of all lines of $\mathbb{R} ^d$. A
line is ``un-sensed'', that is, it is seen as a subset of $\mathbb{R}
^d$, without a preferred direction. For $K$, a compact
of $\mathbb{R} ^d$, the \emph{hitting set of $K$} is the set of lines
that intersect $K$, denoted as
\begin{eqnarray*}
[K] & =& \bigl\{ l \in\mathcal{L} ^d: \mbox{$l$ hits $K$} \bigr\}.
\end{eqnarray*}

We also denote by $m_d$ the Hausdorff measure of dimension $d$. With
this notation, we have the following.
\begin{defin}
\label{lines}
The \emph{invariant line measure} $\mu_d(l)$ is the unique measure on
the space of lines of $\mathbb{R} ^d$ that is invariant under Euclidean
isometries, and normalized by the following requirement: for a compact
set $K$ in $\mathbb{R} ^d$ of nonempty interior, the $\mu_d$-measure
of the hitting set of $K$ is half the Hausdorff $(d-1)$-dimensional
measure of the boundary of $K$:
\begin{eqnarray*}
\mu_d\bigl([K]\bigr) & =& \tfrac{1}{2} m_{d-1}(\partial
K).
\end{eqnarray*}
\end{defin}
The reason for the normalizing constant $\frac{1}2$ is to ensure that
the measure of the hitting set of a flat hypersurface $A$ is its
hyperarea $m_{d-1}(A)$.

We will often need the hitting set of a ball, so write $\omega_{d-1}$
for the hyperarea of the unit sphere $\mathbb{S} ^{d-1}$. Thus $\mu
_{d}([B(x,r)]) = \frac{\omega_{d-1}}{2} r^{d-1}$.\vspace*{1pt}

Consider the following parameterization of a line $l$: it is given by a
direction and a localization. The direction is an element of the
projective space $P \mathbb{R} ^{d-1}$. It then defines a hyperplane
normal to this direction, through a special point---the origin---of
$\mathbb{R} ^d$. The localization is a point on this hyperplane. The
line $l$ is then the line through this point normal to this hyperplane.

Hence, writing $\mathcal{H} = \mathbb{R} ^{d-1}$ for an hyperplane of
$\mathbb{R} ^d$, we may parameterize the set of lines by $P \mathbb{R}
^{d-1} \otimes\mathcal{H} $. Notice that if we want to keep track of
the topology of the set of lines, the product should be twisted, but we
only need measure-theoretical properties, so we stick to the simpler
direct product.

Recall that the projective space $P \mathbb{R} ^{d-1}$ may be seen as
the sphere $\mathbb{S} ^{d-1}$ with opposite---antipodal---points
identified. The projective space then inherits the natural metric on
the sphere, namely the distance between two pairs of antipodal points
is the angle between the pair of segments joining them, in radians. Up
to a null-measure set, the projective space may be more simply viewed
as a hemisphere.

With this parameterization, and writing $\mathbb{B} ^{d}$ for the unit
ball in $\mathbb{R} ^d$, we may write $\mu_d$ as a product measure on
$P \mathbb{R} ^{d-1} \otimes\mathcal{H} $:
\[
\mu_d = \frac{1}{m_{d-1}(\mathbb{B} ^{d-1})} m_{d-1} \otimes m_{d-1}.
\]

To make a clearer reference to it, we write $\mu_{d-1}^{(o)} = 2
m_{d-1} / m_{d-1}(\mathbb{S} ^{d-1})$ for the probability measure on $P
\mathbb{R} ^{d-1}$. For a set of lines $\mathcal{L} $, we write $\mu
_{d-1}^{(o)}( \mathcal{L} )$ for the measure of the set of directions
of lines in $\mathcal{L} $ that go\vspace*{1pt} through the origin. In particular,
for $A$ a subset of $\mathbb{R} ^d$ and $x$ a point of $\mathbb{R} ^d$,
the value $\mu_{d-1}^{(o)}([A - x])$ will be the proportion of
directions that appear in the set of directions of lines through $A$
and $x$. We call it the \emph{solid angle} of $A$ from $x$. Notice that
it is not quite the usual definition since pairs of antipodal points
contribute only once, and since we have normalized to a probability
measure. This is because we are concerned with orientations of
undirected lines, rather than angles of directed lines.

A special case we shall need is the measure of a cone, that is, all the
lines with angle less than $\theta_0$ from a given line through the
origin. This measure is proportional to the area on the hemisphere hit
by the lines, hence proportional to $ \int_{0}^{\theta_0} \sin^{d-2}
\theta \cos\theta\,\mathrm{d}\theta = \sin^{d-1}(\theta_0) / (d-1)$.
The integral up to $\theta_0 = \pi/2$ has value~$1$, so that
%
\begin{eqnarray}
\label{amettreintro} \mu_{d-1}^{(o)} ( \mbox{cone of aperture $\theta
_0$} ) & =& \sin^{d-1}(\theta_0).
\end{eqnarray}

The uniform Poisson line process is the image of the Poisson point
process on~$\mathcal{L} ^d$ with intensity measure $\mu_d$.

We define our improper Poisson line process by adding a mark on each
line, a speed limit. Namely, the improper Poisson line process is the
image of a Poisson point process on $P \mathbb{R} ^{d-1} \otimes
\mathcal{H} \otimes\mathbb{R} ^{+}_*$ with intensity measure $\mu_{d,
\gamma}$, given by the density
\begin{eqnarray*}
\mathrm{d}\mu_{d, \gamma}(l,v) & =& \mathrm{d}\mu_d(l) (\gamma-
1) v^{- \gamma} \,\mathrm{d}v
\end{eqnarray*}
for $\gamma> d$. \citet{Kendall2014} does define this process for
all $\gamma> 1$, but the relevant case for SIRSNs is that of $\gamma> d$.

In words, we have more and more slower lines, following a power law.
Since $\int v^{-\gamma} \,\mathrm{d}v$ diverges at zero, the lines are
dense in $\mathbb{R} ^d$. However, lines faster than any given speed
are not dense. In particular, the number of lines faster than $v_0$
hitting a convex set $K$ is a Poisson variable with parameter
\begin{eqnarray}
\mu_{d, \gamma} \bigl( (l,v): l \in[K]\mbox{ and }v \geq v_0 \bigr)
& =& \mu_d\bigl([K]\bigr) \int_{v_0}^{\infty}
(\gamma- 1) v^{- \gamma} \,\mathrm{d}v
\nonumber
\\
& =& \frac{1}{2} m_{d-1}(\partial K) v_0^{- (\gamma- 1)}.
\nonumber
\end{eqnarray}

We call $\Pi= \Pi(d, \gamma)$ the corresponding random process of
marked lines $(l, v)$. Since the dimension $d$ and parameter $\gamma$
will always be clear from context, we drop them in the notation. Notice
that the total number of lines is almost surely countable. If $(l,v)
\in\Pi$, we say that the speed of line $l$ is $v$ and denote it $v(l)$.

For a subset $\mathcal{L} $ of lines, we write $\Pi_{\mathcal{L} }$
for the restriction of $\Pi$ to these lines, that is, $\Pi_{\mathcal{L}
} = \{ (l,v): l \in\mathcal{L} \} $. In particular, the line process
restricted to lines hitting $A$ but not $B$ is $\Pi_{[A] \setminus[B]}$.

We denote $\mathcal{S} $ the \emph{silhouette} of $\Pi$, that is the
random set in $\mathbb{R} ^d$ made of all the lines of $\Pi$, that is
$\mathcal{S} = \{ x\in\mathbb{R} ^d: \exists(l, v) \in\Pi: x\in l \}
$. We also write $\mathcal{S} _{v_0}$ for the random closed set in
$\mathbb{R} ^d$ made of all the lines $(l,v)$ in $\Pi$ such that $v
\geq v_0$.

\begin{figure}

\includegraphics{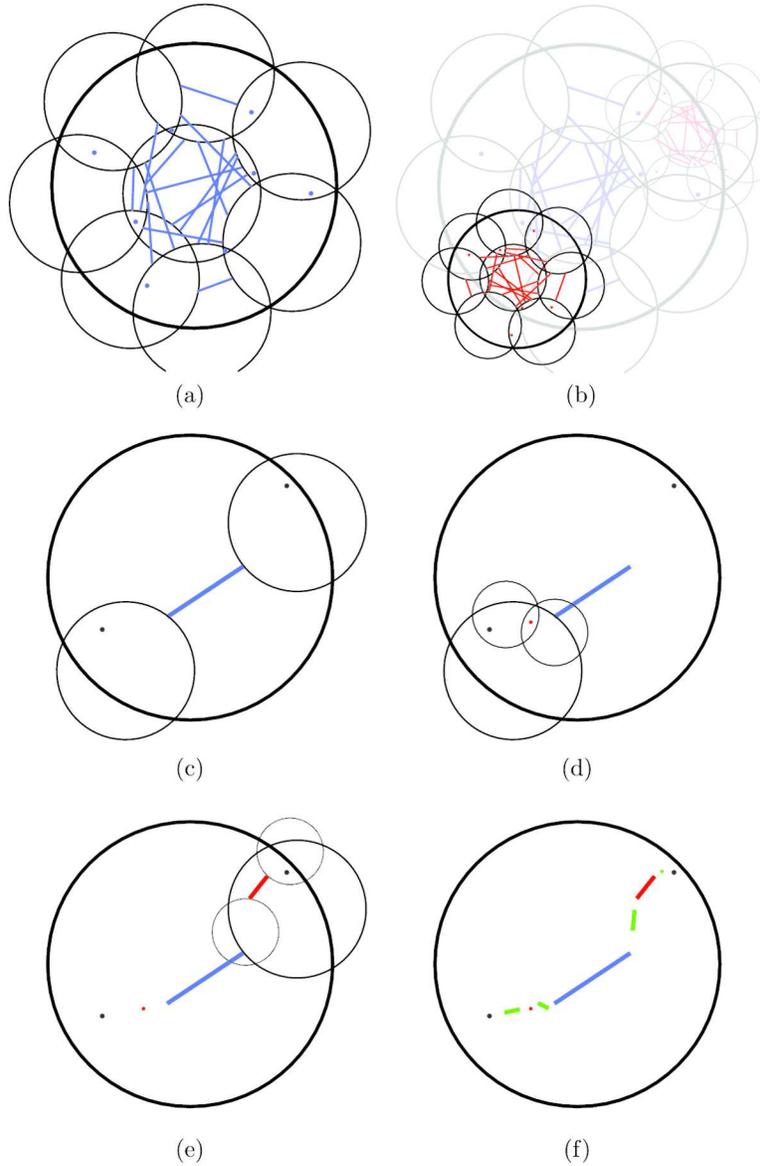}

\caption{We cover each ball $B_i$ with a $r_1$-net, and connect each
pair of balls of the net with a segment, be it a point \textup{(a)}. We
then do the same at each scale, for $r_{n+1}$-nets of the balls of the
$r_n$-nets~\textup{(b)}. Each point belongs to a ball of the $r_1$-net,
so we build a path between two points starting from the segment
connecting their balls~\textup{(c)}. We then \textup{(d)} and \textup{(e)} connect the points
to the segment endpoints with the segment connecting their balls in
their respective $r_2$-nets. At stage $n$, the path \textup{(f)} is
made of $2^n$ segments at each scale $n$.}\label{paths_in_pieces}\label{fleur5}\label{fleur7}
\end{figure}

We may then define $\Pi$-paths.
\begin{defin}
\label{pipath}
A finite-time $\Pi$-path is a locally Lipschitz path in $\mathbb{R}
^d$ respecting the speed limits imposed by $\Pi$. More precisely, it is
a continuous \mbox{$\mathbb{R} ^d$-}valued function
\begin{eqnarray*}
\xi: [0, T] & \to& \mathbb{R} ^d,
\end{eqnarray*}
with $T$ finite, such that for almost all $t \in[0, T]$, either:
\begin{itemize}
\item{the speed is zero: $\xi'(t) = 0$;}
\item{or the path follows a line in $\Pi$: there is a $v \geq\llvert
\xi'(t) \rrvert$ such that $(\xi(t) + \xi'(t) \mathbb{R}, v) \in\Pi$.}
\end{itemize}

We call $T = T(\xi)$ the \emph{time length} of the path $\xi$, or just
its \emph{time} for short.

An infinite-time $\Pi$-path $\xi$ is the same, with $T(\xi) = \infty$,
except that its domain is $[0, \infty)$.

A $\Pi$-path is a finite or infinite-time $\Pi$-path.

In an abuse of notation, we write $\xi\in\Pi$.
\end{defin}
Notice that the image of a $\Pi$-path is not necessarily contained in
the silhouette~$\mathcal{S} $. It only needs to have speed zero outside
$\mathcal{S} $. The remark is especially relevant in dimension at least
three, where the lines never cross. However, since the lines are dense
in $\mathbb{R} ^d$, it turns out that there are paths joining any pair
of points, without any segment
in $\mathbb{R} ^d \setminus\mathcal{S} $. We give a clearer intuition
of their tree-like structure in the proof of Theorem~\ref{Tform} and
Figure~\ref{paths_in_pieces}.

%
%

We write $\Pi^{ab} = \{ \xi\in\Pi: \xi(0) = a$ and $\xi(T(\xi)) = b
\} $ for any two points $a$ and~$b$.

The following theorem is a union of results from \citeauthor
{Kendall2014}'s (\citeyear{Kendall2014}) paper.
%
\begin{teo}
\label{Kendall} Almost surely, all finite-time $\Pi$-paths have finite
Euclidean length.

Almost surely, there are finite-time $\Pi$-paths between each pair of
points of $\mathbb{R} ^d$. Moreover, for any two points $a$ and $b$,
the infimum $T_{ab}$ of time lengths $T(\xi)$ of $\Pi$-paths $\xi\in
\Pi^{ab}$ is attained.

Hence, $\mathbb{R} ^d$ with the metric $d(a,b) = T_{ab}$ is a \emph
{random metric geodesic} space.
\end{teo}
We call this metric \emph{time length} or \emph{$\Pi$-length}. Time
diameters and similar notions are defined in the same way.

We denote by $\mathcal{N} $ the random network made by all the
geodesics connecting all pairs of points in $\mathbb{R} ^d$. Our aim is
to show that $\mathcal{N} $ is a SIRSN.

It is often possible to define similar metrics on other sets $\Gamma$
of marked lines, though they might not be geodesic. We then speak of
$\Gamma$-length. The typical case is when we restrict $\Pi$ to a subset
of lines $\mathcal{L} $, yielding $\Gamma= \Pi_{\mathcal{L} }$ and $\Pi
_{\mathcal{L} }$-length.

We now introduce some notation and remarks to make easier manipulating
paths and geodesics:
\begin{itemize}
\item{If $\xi\in\Pi^{ab}$, we often write $\xi_{ab}$ instead.}
\item{Concatenation of $\Pi$-paths is denoted $\xi_{ac} = \xi_{ab}\xi
_{bc}$, that is, $\xi_{ac}(t) = \xi_{ab}(t)$ if $t < T(\xi_{ab})$, and
$\xi_{ac}(t) = \xi_{bc}(t- T(\xi_{ab}))$ if $t \geq T(\xi_{ab})$. }
\item{We use the letter $g$ for geodesics, and usually $g_{ab}$ for a
geodesic from $a$ to $b$. We say that $g_{ab}$ is unique if there is a
unique geodesic from $a$ to $b$.}
\item{In an abuse of notation, we identify a $\Pi$-path $\xi$ with its
image in $\mathbb{R} ^d$ whenever it is clear. Hence, we may write $x
\in\xi$ if there is a time $t$ such that $\xi(t) = x$. Similarly, if
$g_{ab}$ is unique and $c, d \in g_{ab}$, we say that $g_{cd}$ is
included in $g_{ab}$.}
\item{If $g_{ab}$ is unique and $c \in g_{ab}$, then $g_{ab} = g_{ac} g_{cb}$.}
\item{$T(g_{ab}) = T_{ab}$.}
\item For a line $l \in\Pi$ and a $\Pi$-path $\xi: [0, T] \to\mathbb
{R} ^d$, we define the \emph{intersection length} of $\xi$ and $l$ as
$L_{\xi}(l) \,\hat{=}\, m_1(l \cap\xi([0,T])) $.
\item{If the intersection length of $\xi$ and $l$ is not zero, we say
that $l$ is in the \emph{support} of $\xi$, or that it is supporting
$\xi$. We write $l \sqsubset\xi$. Moreover, we denote the support of
$\xi$ by $\mathcal{L} _{\xi} \,\hat{=}\, \{ l \in\Pi: L_{\xi}(l) > 0 \} $.}
\item{We define the \emph{intersection time} of $l$ and $\xi$ as
$T_{\xi}(l) \,\hat{=}\, m_1([0,T] \cap\xi^{-1}(l)) $. }
\item{In particular, almost surely, for all geodesics $g$, we have the
following equality and decompositions:
%
\begin{eqnarray}\label{ll}
L_{g}(l) & =& v(l) T_g(l) \qquad\mbox{for
all }l \in\Pi,
\nonumber
\\
T(g) & =& \sum_{l \in\mathcal{L} _{g}} T_g(l) 
\\
& =& \sum_{l \in\mathcal{L} _{g}} v(l) L_g(l).\nonumber
\end{eqnarray}
}
\item{If the support of $\xi$ is included in $\mathcal{L} $, that is
$\mathcal{L} _{\xi} \subset\mathcal{L} $, we say that $\xi$ is $\Pi
_{\mathcal{L} }$-path. We abuse notation by writing $\xi\in\Pi
_{\mathcal{L} }$.}
\item{Similarly, we write $\xi\in\Pi_{\mathcal{L} }^{ab}$ if $\xi\in
\Pi_{\mathcal{L} }$ and $\xi\in\Pi^{ab}$.}
\end{itemize}

Finally, a few more generic notation. We call \emph{internal
$\varepsilon$-net} of a subset $A$ of $\mathbb{R} ^d$ any maximal
subset $x_1, \dots, x_k$ of $A$ such that $\llvert x_i - x_j \rrvert
\geq\varepsilon$ for all $i \neq j$. For a set $A \subset\mathbb{R}
^d$, we write $A^r$ for its $r$-widening, that is the Minkowski sum $
A^r = A \oplus B(0, r)$. We denote the maximal speed of a set of lines
$\mathcal{L} $ by $v_{\max}(\mathcal{L} )$. That is, $v_{\max}(\mathcal
{L} ) = \sup_{l \in\mathcal{L} } v(l)$. Notice that if those lines all
hit a compact set $K$, this supremum is a maximum. We abuse notation by
writing $v_{\max}(A) = v_{\max}([A])$ for $A \in\mathbb{R} ^d$. We use
$C, c, c_1, \dots, c_i$ for any positive constant.

\section{\texorpdfstring{$\Pi$}{$Pi$}-diameters of sets}
\label{time_diam}

We start with giving a more quantitative version of Theorem $3.6$ in
\citeauthor{Kendall2014}'s (\citeyear{Kendall2014}) article. Namely, we
show that in a given precompact set $A$, any two points can be joined
in finite time, and that the largest time between two such points---the $\Pi$-diameter of $A$---is not too big with high probability: this
random variable has more than an exponential moment. We include a
generalization that we will need later on, by allowing the possibility
of ignoring lines hitting forbidden areas $F$.

%
\begin{teo}
\label{Tform}
Recall that $d \geq2$ and $\gamma> d$.
Let $\alpha_{\min} = 2^{(\gamma- 1) / (\gamma- d)}$. Let $\Omega<
\Omega_{\max} = (4 \alpha_{\min} )^{1 - d} $ and define $\alpha_{\max} =
\Omega^{-\afrac{1}{d-1}} / 4$. Choose $\alpha$ such that $\alpha
_{\min} < \alpha < \alpha_{\max}$, and note in particular that $\alpha
> 1$.

Let $A$ and $F$ be two subsets of $\mathbb{R} ^d$, such that, for some
$r > 0$:
\begin{itemize}
\item{$A$ is connected.}
\item{From any $x\in A^{r / (\alpha-1) }$, the solid angle $\mu
^{(o)}_{d-1}([F - x])$ of $F$ is less than $\Omega$.}
\item{$A$ may be covered by $\mathring{N}$ balls $B_i$ of radius $r$.}
\end{itemize}

Then there is a $T_1$ depending only on $\alpha$, $\Omega$, $\gamma$
and $d$ such that, for any $ \varepsilon_{\max} = 1 / (2 \mathring{N} (2
\alpha+1)^d) > \varepsilon> 0$, with probability $1 - \varepsilon/
\varepsilon_{\max} $, the diameter of $A$ on $\Pi_{[A^{r / (\alpha-
1)}] \setminus[F]}$ satisfies the bound
%
\begin{eqnarray}
T_{A, F} & \,\hat{=}& \sup_{x,y \in A} \inf
_{\xi\in\Pi^{xy}_{[A^{r / (\alpha- 1)}]\setminus[F]}} T(\xi)
\nonumber
\\
\label{td} & \leq& T_{r, \mathring{N}} \biggl( \ln\frac{1 }{\varepsilon}
\biggr)^{\afrac{1}{\gamma-1} } \qquad\mbox{with}
\\
T_{r, \mathring{N}} & =& \mathring{N} T_1 r^{\vafrac{\gamma- d}{\gamma
-1}}.\nonumber
\end{eqnarray}
In particular, this maximal time has all exponential moments, and more:
for any $\delta< T_{r, \mathring{N}}^{-\afrac{1}{\gamma- 1}} $, we have
%
\begin{eqnarray}
\label{morexp} \mathbb{E} \bigl[ \exp\bigl( \delta T_{A,F}^{\gamma- 1 }
\bigr) \bigr] & <& \infty.
\end{eqnarray}
\end{teo}

The proof is a slight variation on that of Theorem $3.6$ in \citeauthor
{Kendall2014}'s (\citeyear{Kendall2014}) article.

\begin{pf*}{Proof of Theorem \ref{Tform}}
Since $A$ is connected and covered by $\mathring{N}$ open balls $B_i$
of radius $r$, we may build a path between any two points of $A$ by
concatenating at most $\mathring{N}$ paths between two points of $A$
belonging to the same ball $B_i$ of the cover.

We now recursively build a path between each pair of points $x_0$ and
$y_0$ of $B_i$, in a binary tree-like fashion. First, we specify
\begin{eqnarray*}
r_n & =& r \alpha^{-n}.
\end{eqnarray*}

We will choose corresponding speed limits $v_n$ later. Given such $v_n$:
\begin{itemize}
\item{ We call $A_0$ the set of balls $\{ B_i\}_{1\leq i \leq\mathring{N}}$.}
\item{To any ball $B^{(n)} \in A_n$, we associate an internal
$r_{n+1}$-net of that ball. It may be viewed as a collection of balls
$B^{(n+1)}$ of radius $r_{n+1}$. We then define $A_{n+1}$ as the set of
all these balls $B^{(n+1)}$ for all $B^{(n)}$ together.}
\item{We have thus built nested internal $r_n$-nets.}
\item{For any two balls $B_i^{(n+1)}$ and $B_j^{(n+1)}$ belonging to
the internal $r_{n+1}$-net of the same ball $B^{(n)}$ in $A_{n}$, we
find a line of speed at least $v_{n}$ that hits both $B_i^{(n+1)}$ and
$B_j^{(n+1)}$, but\vspace*{1pt} not $F$. We will have to prove this is possible with
high enough probability.}
\item{For $x_{n}$ and $y_{n}$, both belonging to $B^{(n)} \in A_{n}$,
we may then find two points $x_{n+1}$ and $y_{n+1}$ such that: $x_{n}$
and $x_{n+1}$ (resp., $y_{n}$ and $y_{n+1}$) belong to the same ball
$B_i^{(n+1)} \in A_n$ (resp., $B_j^{(n+1)}$), and $x_{n+1}$ and
$y_{n+1}$ both belong to the same line of speed at least $v_{n}$.
}
\item{We may then build a path between $x_{n}$ and $y_{n}$ as a
concatenation of three paths: $\xi_{x_n y_n} = \xi_{x_n x_{n+1}} \xi
_{x_{n+1} y_{n+1}} \xi_{y_{n+1} y_n} $. The middle one is a segment.
The other two are paths between points of the same ball in $A_{n+1}$.}
\end{itemize}

As illustrated in Figure~\ref{paths_in_pieces}, we thus obtain a path
between $x_0$ and $y_0$ that is made of exactly $2^n$ segments for each
$n$, each at speed at least $v_{n}$ between two balls of the same
internal $r_{n+1}$-net
of a ball of radius $r_{n}$. Moreover, since the points $x_0$ and $y_0$
are in $A$, all of the segments are between points of the Minkowski sum
$ A \bigoplus_{n=1}^{\infty} B(0, r_n) = A^{r / (\alpha- 1)}$.

This construction has built a path for each pair of points $x_0$ and
$y_0$ in $B_i$. Since segments between balls of the same internal $r_{n+1}$-net
of a ball of radius $r_{n}$ are at most $(2 r_n + 2 r_{n+1})$ long, and
are at speed at least $v_n$, the $\Pi_{ [A^{r / (\alpha- 1)}]
\setminus[F]}$-diameter $T_{A,F}$ of $A$ is bounded from above by
%
\begin{eqnarray}
\label{gen_bound_A} T_{A, F} & \leq&\mathring{N} \sum
_{n=1}^\infty2^n \frac{2 r_n + 2 r_{n+1} }{v_n},
\end{eqnarray}
on the event that this construction is possible.

Now to control the probability of this event, we need:
\begin{itemize}
\item{a bound on the number of pairs of balls in the same $r_n$-net;}
\item{a bound on the probability that the fastest line hitting two such
balls but not $F$ is slower than $v_n$.}
\end{itemize}

We obtain the first bound by using the formula $r_n = \alpha^{-1} r_{n-1}$.
Indeed, each $r_n$-net of a $r_{n-1}$-ball is then the same as a
$\alpha^{-1}$-net of a radius $1$ ball. Since the balls $B(s_i, 1/(2
\alpha))$ centered on the points of such a $\alpha^{-1}$-net are
disjoint, and all included in a ball of radius $1 + 1/(2 \alpha)$, a
volume argument shows that there at most $(2 \alpha+ 1)^d$ balls in
each $r_n$-net. So that there are at most $\mathring{N}(2\alpha
+1)^{dn}$ balls in $A_n$, and at most $\mathring{N} (2\alpha+1
)^{d(n+1)}$ pairs of balls in the same $r_n$-net.

We now consider these lines that hit both of two $r_{n+1}$-balls $B(x,
r_{n+1})$ and $B(y, r_{n+1})$ in an internal $r_{n+1}$-net of a
$r_n$-ball, but that do not hit $F$. We have seen that they were in $
[A^{r / (\alpha- 1)} ]$. We may then use the hypothesis of the theorem
on the solid angle of $F$.

We want a bound on
\[
\mu_d \bigl( \bigl[ B(x, r_{n+1}) \bigr] \cap\bigl[ B(y,
r_{n+1}) \bigr] \cap[ F ]^c \bigr).
\]

%
\begin{figure}

\includegraphics{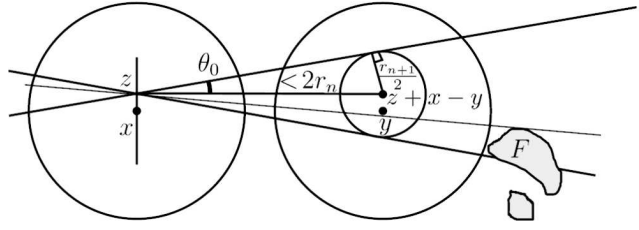}

\caption{Lines that hit both balls $B(x, r_{n+1})$ and $B(y, r_{n+1})$
but not $F$ (not necessarily connected) contain lines through $z \in
D(x, r_{n+1}/2)$ and $B(z + x - y, r_{n+1} / 2)$. These are included in
the cone of angle $\theta_0$ through $z$.}
\label{disques}
\end{figure}

Now, as illustrated in Figure~\ref{disques}, let $D(x, r_{n+1} / 2)$ be
the hyperdisk centered on $x$ with radius $r_{n+1} / 2$ and included in
the hyperplane perpendicular to the vector $x - y$. It is included in
$B(x, r_{n+1} )$. Moreover, for any $z \in D(x, r_{n+1} / 2)$, the ball
$B(z + y - x, r_{n+1} / 2)$ is included in $B(y, r_{n+1})$.

Now we may write
\begin{eqnarray*}
&& \mu_d \bigl( \bigl[ B(x, r_{n+1}) \bigr]
\cap\bigl[ B(y, r_{n+1}) \bigr] \cap[ F ]^c \bigr)
\\
&&\qquad  \geq m_{d-1} \bigl( D(x, r_{n+1} / 2) \bigr)
\\
&&\quad\qquad{}\times\inf
_{z \in D(x, r_{n+1} / 2) } \mu_{d-1}^{(o)} \bigl( \bigl[ B( y -
x, r_{n+1} / 2) \bigr] \cap[F-z]^c \bigr)
\\
&&\qquad  \geq m_{d-1} \bigl( D(x, r_{n+1} / 2) \bigr) \times\bigl(
\mu_{d-1}^{(o)} \bigl( \bigl[ B(y - x, r_{n+1} / 2)
\bigr] \bigr) - \Omega\bigr).
\end{eqnarray*}

Now, $y - x \leq2 r_n = 2 \alpha r_{n+1}$. So that the lines in $ [
B(y - x, r_{n+1} / 2) ] $ are in a cone of angle $\theta_0$ such that
$ \sin(\theta_0) \leq(4 \alpha)^{-1}$. So that by formula~(\ref{amettreintro}):
%
\begin{eqnarray}
\label{aperture} \mu_{d-1}^{(o)} \bigl( \bigl[ B(y - x,
r_{n+1} / 2) \bigr] \bigr) & \leq&( 4 \alpha) ^{1 - d}.
\end{eqnarray}

Since $\alpha< \alpha_{\max}$, we have $ ( 4 \alpha ) ^{1 - d} -
\Omega > 0$.

Moreover, with $\nu_{d-1}$ the $(d-1)$-volume of a unit $(d-1)$-ball,
we have
\begin{eqnarray*}
m_{d-1} \bigl( D(x, r_{n+1} / 2) \bigr) & =& \nu
_{d-1} \biggl( \frac{r_{n+1}}{2 } \biggr) ^{d-1},
\\
\mu_d \bigl( \bigl[ B(x, r_{n+1}) \bigr] \cap\bigl[ B(y,
r_{n+1}) \bigr] \cap[ F ]^c \bigr) & \leq&
\nu_{d-1} \bigl( ( 4 \alpha) ^{1 - d} - \Omega\bigr) \biggl(
\frac{r_{n+1}}{2 } \biggr) ^{d-1}
\\
& \,\hat{=}&  c_{\alpha, \Omega } r_{n+1}^{d-1},
\end{eqnarray*}
using the notation $c_{\alpha, \Omega }$ to summarize and emphasize
that this bound depends on $\alpha$ and $\Omega$, but not on $n$,
$r$, or any other feature of $A$ and $F$.

We thus get
\begin{eqnarray*}
\mathbb{P} \bigl[v_{\max} \bigl( \bigl[ B(x, r_{n+1}) \bigr]
\cap\bigl[ B(y, r_{n+1}) \bigr] \cap[ F ]^c \bigr) \leq
v_n \bigr] & \leq&\exp\biggl( - c_{\alpha, \Omega } \frac{
r_{n+1}^{d-1}}{v_n^{\gamma-1}}
\biggr).
\end{eqnarray*}

Multiplying by the number of relevant pairs of ball, we then obtain
that the construction is possible except on an event of probability at most
%
\begin{equation}
\label{prob_cons_bound} \sum_{n=0}^{\infty}
\mathring{N} (2 \alpha+1)^{d(n+1)} \exp\biggl( - c_{\alpha, \Omega }
\frac{ r_{n+1}^{d-1}}{v_n^{\gamma-1}} \biggr).
\end{equation}

Taking
\begin{eqnarray*}
v_n & =& \frac{r_{n+1}^{(d-1) / (\gamma- 1)}}{ ( (n+1) c_{\alpha,
\Omega} \ln(1/\varepsilon) ) ^{1/(\gamma- 1)}}
\end{eqnarray*}
and replacing in bound \eqref{prob_cons_bound}, we see that this
becomes less than
\[
\sum_{n=0}^{\infty} \mathring{N} (2 \alpha
+1)^{d(n+1)} \varepsilon^{n+1} \leq\frac{\varepsilon}{\varepsilon_{\max}},
\]
for $\varepsilon< 1 / (2 \mathring{N} (2 \alpha+1)^{d(n+1)} ) =
\varepsilon_{\max}$, so that $\varepsilon_{\max}$ does not depend on
$A$ or $F$ except through $\alpha$ and $\mathring{N}$.

Replacing $v_n$ and $r_n$ with their value in the bound \eqref
{gen_bound_A}, we get
\begin{eqnarray*}
T_{A,F} & \leq&\mathring{N} \sum_{n=1}^\infty
2^{n+1} (1 + \alpha) \frac{r_{n+1} }{v_n}
\\
& =& \mathring{N} \Biggl( (1 + \alpha) c_{\alpha, \Omega}^{\afrac
{1}{\gamma-1} } \sum
_{n=2}^\infty\bigl(2 \alpha^{-\vafrac{\gamma-d}{\gamma-1} }
\bigr)^{n} n^{\afrac{1}{\gamma-1}} \Biggr)
\\
&&{}\times r^{(\gamma- d) / (\gamma-
1)} \biggl( \ln
\biggl( \frac{1}{\varepsilon} \biggr) \biggr) ^{\afrac{1}{\gamma-1} },
\end{eqnarray*}
where the first parentheses correspond to $T_1$ and $T_1$ is finite
since $\alpha> \alpha_{\min}$.
We have thus proved formula~\eqref{td} of the theorem.

The moment \eqref{morexp} is simple integration:
\begin{eqnarray*}
\mathbb{E} \bigl[ \exp\bigl( \delta T_{A,F}^{\gamma- 1 } \bigr)
\bigr] & \leq& C + \frac{1}{\varepsilon_{\max}} \int_{0}^{\varepsilon_{\max}}
\exp\biggl( \delta\biggl(T_{r, \mathring{N}} \biggl( \ln\frac
{1}{\varepsilon}
\biggr) ^{\afrac{1}{\gamma- 1}} \biggr)^{\gamma- 1 } \biggr) \,\mathrm
{d}\varepsilon
\\
& \leq& C + \frac{1}{\varepsilon_{\max}} \int_{0}^{\varepsilon_{\max}}
\varepsilon^{- \delta T_{r, \mathring{N}}^{\afrac{1}{\gamma- 1}}}\,
\mathrm{d}\varepsilon
\\
& <& \infty.
\end{eqnarray*}
Here, we use $\delta< T_{r, \mathring{N}}^{-\afrac{1}{\gamma- 1}}$ and
$ C $ corresponds to the integral between $\varepsilon_{\max}$ and~$1$,
bounded by the value of $T_{A,F}$ for $\varepsilon= \varepsilon_{\max}$.
\end{pf*}

%
\begin{rem}\label{...}
If no lines are forbidden, that is if $F = \varnothing$, then
$\alpha$ can be taken as big as we wish, so that the lines used all
hit as small a widening of $A$ as we want.

We may slightly generalize the theorem: instead of using a
forbidden area $F$, we could use different conditions for which lines
to accept. The important property is that we must have enough relevant
lines hitting pairs of balls in a $r_{n+1}$-net.

There are a few optimisations that could be used to gain slightly
in the constants. For example, we have written the proof with one
subnet inside each ball. If we had used a single $r_n$-net of a
correctly widened $A$, we would have only about $\mathring{N} 2^d
\alpha^{d(n+1)}$ pairs of ball, allowing a bigger $\varepsilon
_{\max}$. The result stays essentially the same, however, as proven by
the following proposition.
\end{rem}

%
\begin{prop}
\label{optimal_diameter}
For any two points $x$ and $y$, their $\Pi$-distance $T_{xy}$ does not
have a moment with higher exponent on the time $T_{xy}$ than in formula~(\ref{morexp}). That is, for any $\delta> 0$, for any $\eta> \gamma
- 1$, we have
\begin{eqnarray*}
\mathbb{E} \bigl[ \exp\bigl( \delta T_{xy}^{\eta} \bigr)
\bigr] & =& \infty.
\end{eqnarray*}
\end{prop}

\begin{pf}
Say that $\llvert x - y \rrvert = r$. A path from $x$ to $y$ has to
go from $x$ to the border of $B(x,r)$. Hence, it must use lines hitting
$B(x,r)$ for a length at least $r$. So that $T_{xy}$ is controlled by
the fastest line hitting $B(x,r)$. Now,
\begin{eqnarray*}
\mathbb{P} \bigl[ v_{\max}\bigl(B(x, r)\bigr) \leq v_{\varepsilon}
\bigr] & =& \exp\bigl( - c(r) v_{\varepsilon}^{1 - \gamma} \bigr)
 = \varepsilon
\end{eqnarray*}
with
\begin{eqnarray*}
v_{\varepsilon} & =& \biggl( \frac{1}{c(r)} \ln\frac{1}{\varepsilon}
\biggr)
^{\afrac{1}{1 - \gamma} },
\end{eqnarray*}
where $c(r)$ depends only on $r$.

So that on an event of probability at least $\varepsilon$, we have
the bound $T_{xy} \geq r / v_{\varepsilon} = c_2(r) (\ln(1 /
\varepsilon))^{1 / (\gamma- 1)}$. Hence, for some positive constant
$c_3$ depending on $\delta$, $\eta$ and~$r$:
\begin{eqnarray*}
\mathbb{E} \bigl[ \exp\bigl( \delta T_{xy}^{\eta} \bigr)
\bigr] & \geq&\int_{0}^1 \exp\biggl( -
c_3 \biggl( \ln\frac{1}{\varepsilon} \biggr)^{\eta/ (\gamma- 1)}
\biggr)\,\mathrm{d} \varepsilon
\\
& =& \infty,
\end{eqnarray*}
since $\eta/ (\gamma- 1) > 1$.
\end{pf}

\section{Almost sure uniqueness of \texorpdfstring{$\Pi$}{$Pi$}-geodesics}
\label{unique}

In this section, we prove that the network $\mathcal{N}$ satisfies the
property~1 of a SIRSN, that is that between two specified
points $x$ and $y$ in $\mathbb{R} ^d$, there is almost surely only a
unique route. Since routes are the $\Pi$-geodesics, this is equivalent
to almost sure uniqueness of the geodesic $g_{xy}$.

The case in dimension $2$ has already been established in \citeauthor
{Kendall2014}'s (\citeyear{Kendall2014}) paper. The following proof, on
the other hand works in all dimensions more than $2$, as stated in
Theorem~\ref{unique_geodesics}. Note that it does not work in dimension $2$.

The strategy is the following:
\begin{itemize}
\item{We introduce a concept of \emph{many directions}, with the
following property. If a $\Pi$-path has many directions near a point
$x$, any finite $\Pi$-path supported by the same lines contains $x$.}
\item{We show that almost surely all geodesics have many directions
near all the ends of their component segments, except for the two
extremal points. This is the step where $d\geq3$ is needed.}
\item{We show that for specified $x$ and $y$ in $\mathbb{R} ^d$, almost
surely all geodesics from $x$ to $y$ are supported exactly on the same lines.}
\item{So that almost surely all such geodesics contain the same segment
ends, and this will prove they are the same.}
\end{itemize}

The author thinks the proof is very technical for something that looks
clear enough, but could not find an easier way. Maybe the need to work
with tree-like paths in dimension at least three is the reason why
there is no obvious argument. Hopefully, the concept of many directions
can be useful elsewhere.

We first state two technical lemmas we need for the proofs in
Section~\ref{tech}. We introduce the notion of ``many directions'' and
give some cases where paths have many directions in Section~\ref
{manydir}, culminating in Lemma~\ref{tour}. We prove that geodesics
must use the same lines in Section~\ref{samelines}, and end the proof
of uniqueness in Section~\ref{asunique}.

\subsection{Technical lemmas}
\label{tech}

The first lemma yields a control on the proportion of balls in a nested
set that are hit by lines faster than a threshold appropriately scaling
with their size.
%
\begin{lem}
\label{nested_balls}
Let $d \geq2$.
Let $\alpha> 1$ be a scale factor. We write $p = \alpha^{1 - d} <
1$. Let $B_i = B(x_i, r_i)$ for $1 \leq i \leq n$ be a set of nested
balls with $r_i = r_0 / \alpha^{i}$ for some $r_0$. For some $v_0$, define
\begin{eqnarray*}
v_i & =& v_0 p^{\afrac{i}{\gamma-1} },
\\
\mathcal{V} _0 & =& \bigl\{ (l, v) \in\Pi: v \geq
v_0 \bigr\},
\\
\mathcal{V} _i & =& \bigl\{ (l, v) \in\Pi: v_{i-1} \geq v
\geq v_i \bigr\}.
\end{eqnarray*}
Define now $H_i$ as the number of balls $B_{i+k}$ smaller than $B_i$
that are hit by lines faster than $v_i$, and $I_n$ as the number of
balls $B_i$ that are hit by lines faster than $v_{i-1}$, that is,
\begin{eqnarray*}
H_i & =& \# \bigl\{ k > 0: \mathcal{V} _i \cap
[B_{i+k}] \neq\varnothing\bigr\},
\\
I_n & =& \# \bigl\{ i \leq n: \exists k > 0: [B_i] \cap
\mathcal{V} _{i - k} \neq\varnothing\bigr\} \leq\sum
_{i=0}^{n-1} H_i.
\end{eqnarray*}
Then $\lambda\,\hat{=}\, \mathbb{E} [ \#\{\mathcal{V} _i \cap[B_i]\} ] =
(1 - p) \omega_{d-1} r_0^{d-1} / (2 v_0^{\gamma- 1})$ independently
of $i\geq1$. Then, for any $\delta> 0$:
%
\begin{eqnarray}
\label{res_tower} \mathbb{P} [ I_n > \delta n ] & \leq&\biggl( ( 1 +
\lambda\delta) \biggl( \frac{1 + \delta}{\delta} \biggr) ^{\delta}
p^{\delta} \biggr) ^n.
\end{eqnarray}
Precise but more cumbersome bounds are given in equations \eqref
{Markov} and \eqref{precis}.
\end{lem}

\begin{pf}
We have
\begin{eqnarray*}
\mathbb{E} \bigl[ \#\bigl\{\mathcal{V} _i \cap[B_{i+k}]
\bigr\} \bigr] & =& \mathbb{E} \bigl[ \#\bigl\{\mathcal{V} _i \cap
[B_{i}]\bigr\} \bigr] p^k
\\
& =& \lambda p^k \qquad\mbox{for all }i \geq1\mbox{ and }k\geq0,
\\
\mathbb{E} \bigl[ \#\bigl\{\mathcal{V} _0 \cap[B_{k}]
\bigr\} \bigr] & =& \lambda p^k \sum_{m=0}^{\infty}
p^m
\\
& =& \frac{1}{1-p} \lambda p^k \qquad\mbox{for all }k \geq0.
\end{eqnarray*}

Thus,
\begin{eqnarray*}
\mathbb{P} [ H_i = 0 ] & \leq&1  \qquad\mbox{for all }i \geq0,
\\
\mathbb{P} [ H_i = k > 0 ] & =& \mathbb{P} \bigl[ \mathcal{V}
_i \cap[ B_{i + k} ] \neq\varnothing\mbox{ and }
\mathcal{V} _i \cap[ B_{i + k + 1} ] = \varnothing\bigr]
\\
& =& \exp\bigl( - \lambda p^{k+1} \bigr) - \exp\bigl( - \lambda
p^k \bigr)  \qquad\mbox{using }[B_{i+1}]
\subset[B_i]
\\
& \leq&1 - \exp\bigl( - \lambda\bigl(p^k - p^{k+1}\bigr)
\bigr)
\\
& \leq&\lambda p^k (1 - p)  \qquad\mbox{for }i \geq1.
\end{eqnarray*}
Similarly,
\begin{eqnarray*}
\mathbb{P} [ H_0 = k > 0 ] & =& \exp\biggl( - \frac{\lambda p^{k+1}}{1
- p}
\biggr) - \exp\biggl( - \frac{\lambda p^k}{1 - p} \biggr)
\\
& \leq&\lambda p^k.
\end{eqnarray*}

Hence, for all $0 < a < 1/p$:
%
\begin{eqnarray}
\mathbb{E} \bigl[ a^{H_i} \bigr] & \leq&1 + \sum
_{k=1}^{\infty} \lambda p^k (1 - p)
a^k
\nonumber
\\
& =& 1 + \lambda(1 - p) \frac{ap}{1 - ap}  \qquad\mbox{for all }i \geq1,
\nonumber
\\
\mathbb{E} \bigl[a^{H_0} \bigr] & \leq&1 + \sum
_{k=1}^{\infty} \lambda p^k a^k
\nonumber\\[-8pt]\\[-8pt]\nonumber
& =& 1 + \lambda\frac{ap}{1 - ap},
\\
\mathbb{E} \bigl[ a^{I_n} \bigr] & \leq&\mathbb{E} \bigl[
a^{\sum_{i=0}^{n-1} H_i} \bigr] \vee1
\nonumber
\\
& \leq&\biggl( 1 + \lambda(1 - p) \frac{ap}{1 - ap} \biggr) ^n
\frac{1 + ap (\lambda- 1)}{1 + ap (\lambda(1 - p) - 1) }  \label{boundpre}
\\
& \leq&\biggl( 1 + \lambda\frac{ap}{1 - ap} \biggr)^n.
\nonumber
\end{eqnarray}

We now use the Markov inequality:
%
\begin{eqnarray}
\label{Markov} \mathbb{P} [ I_n > \delta n ] & \leq&
\frac{\mathbb{E} [a^{I_n} ]}{a^{ \delta n}},
\end{eqnarray}
and optimize upon $a$. Exact optimization requires solving a degree-two
equation and yields a cumbersome solution, so we shall only use here
the solution for the second approximation when $\lambda = 1$, that is
$a = \delta/ ((1 + \delta) p )$, so that \mbox{$ap / (1 - ap) = \delta$}.
Using \eqref{boundpre}, we get
%
\begin{eqnarray}
\mathbb{P} [ I_n > \delta n ] & \leq&\frac{1 + \lambda \afrac{ap}{1 -
ap}}{ 1 + \lambda(1 - p) \afrac{ap}{1 - ap} } \biggl( 1 +
\lambda(1 - p) \frac{ap}{1 - ap} \biggr) ^n a^{- \delta n}
\nonumber
\\
\label{precis} & =& \frac{1 + \lambda\delta}{1 + \lambda(1 - p)
\delta} \biggl( \bigl( 1 + \lambda(1 - p) \delta
\bigr) \biggl( \frac{1 + \delta}{\delta} \biggr) ^{\delta} p^{\delta
} \biggr)
^n
\\
& \leq&\biggl( ( 1 + \lambda\delta) \biggl( \frac{1 + \delta}{\delta
} \biggr)
^{\delta} p^{\delta} \biggr) ^n.
\nonumber
\end{eqnarray}\upqed
\end{pf}

The second lemma gives a guarantee that except on exceedingly rare
events; no significant fraction of uniformly randomly oriented lines
are clumped together in a small number of directions.

%
\begin{lem}
\label{angles}
Let $(\mathcal{X}, d, \mu)$ be a metric space with a probability
measure $\mu$, and such that for all positive and small enough
$\varepsilon$, for any point $x$ in $\mathcal{X} $, the volume of the
ball is bounded in this way:
%
\begin{equation}
\label{double} c_1 \varepsilon^{d-1} \leq\mu\bigl[ B(x,
\varepsilon) \bigr] \leq c_2 \varepsilon^{d-1},
\end{equation}
for $c_1$ and $c_2$ constants depending only on the space $(\mathcal
{X}, d, \mu)$.

Let $\alpha> \beta> \delta> 0$. Consider $ \{ s_i \} _{i\in I}$ a $
n^{-\eta}$-net of $\mathcal{X} $. Consider $ \{ x_j \}_{j\leq\alpha
n} $ $\alpha n$ random $\mu$-i.i.d. points in $\mathcal{X} $.

Then there is no subset of $\beta n$ points in $ \{ x_j \} $ that are
all contained in at most $\delta n$ balls of the net $B(s_i, n^{-\eta
})$, except on an event of sub-exponential probability, at most $ ( C
n^{(\beta- \delta) (\eta(d-1) - 1)} )^{-n} $ with $C$ depending only
on $\alpha, \beta$ and $\delta$, $c_1$ and $c_2$.\vspace*{1pt}

Moreover, the projective space $(P \mathbb{R} ^{d-1}, \theta, \mu
^{(o)}_{d-1})$ satisfies the hypotheses for $(\mathcal{X}, d, \mu)$.
Here $\theta(l_1, l_2)$ is the angle between two lines, and $\mu
^{(o)}_{d-1}$ is the natural probability measure on $P \mathbb{R}
^{d-1}$, defined in Section~\ref{notations}.
\end{lem}

\begin{pf}
Since all the points in the $n^{- \eta}$-net are in disjoint balls of
radius $n^{-\eta}/2$, there are at most $c_3 n^{\eta(d-1) }$ points in
the net, for $n$ big enough.
Now, with $c_4, c_5, c_6$ depending on $\alpha, \beta$ and $\delta$:
\begin{eqnarray*}
&& \mathbb{P} \bigl[\mbox{There are }\beta n\mbox{ points
}x_j\mbox{ all contained in }\delta n\mbox{ balls }B
\bigl(s_i, n^{-\eta}\bigr)\mbox{ of the net} \bigr]
\\
&&\qquad  \leq\# \{ \mbox{subsets of $\beta n$ points} \}\, \#\, \{ \mbox{subsets
of $
\delta n$ balls} \}
\\
&&\qquad\quad{}\times \mathbb{P} [ \mbox{$\beta n$ specific points are all
contained in $\delta n$ specific balls} ]
\\
&&\qquad \leq\pmatrix{\alpha n \cr\beta n} \pmatrix{c_3 n^{\eta(d-1)} \cr
\delta n} \bigl( \delta n c_2 n^{- \eta(d-1)} \bigr)^{\beta n}
\\
&&\qquad  \leq 2^{\alpha n} \bigl(c_4 n^{\eta(d-1) - 1}
\bigr)^{\delta n} \bigl(c_5 n^{1 - \eta(d-1)}
\bigr)^{\beta n}
\\
&&\qquad = \bigl( c_6 n^{(\beta- \delta) (\eta(d-1) - 1)} \bigr)^{-n}.
\end{eqnarray*}

In the projective space, where points are seen as lines through the
origin in~$\mathbb{R} ^d$, a ball of radius $\theta_0$ is exactly a
cone of aperture $\theta_0$. We recall formula~(\ref{amettreintro})
and $ \theta/ 2 \leq\sin(\theta) \leq\theta$ for positive small
$\theta$, and we conclude that property~\eqref{double} is indeed satisfied.
\end{pf}

\subsection{Many directions}
\label{manydir}

Having many directions near a point $x$ intuitively means that the
lines used near the point have so many different unit vectors that the
only way to touch all those lines (a \emph{tour}) with a finite curve
is by touching most of them near $x$.

\begin{defin}
\label{deftour}
For a set of lines $\mathcal{L} = \{l_j\}_{j\in J}$ and a subset of
the Euclidean space $X \subset\mathbb{R} ^d$, a \emph{$\mathcal{L}
$-tour in $X$} is a curve $f$ in $X$ such that for all $j\in J$, there
is a $t_j$ such that $f(t_j) \in l_j$. If $f$ is rectifiable, the tour
is said to be finite; else it is infinite.
\end{defin}

Recalling the notation $\mathcal{L} _{\xi}$ for the support of $\xi$:
\begin{defin}
\label{many_dir_label}
A finite $\Pi$-path $\xi$ has \emph{many directions} near a point $x$
if, for all $\varepsilon> 0$, all $\mathcal{L} _{\xi}$-tours in
$\mathbb{R} ^d \setminus B(x, \varepsilon)$ are infinite.
\end{defin}

As a remark, this concept is only interesting in dimension at least
three: in dimension two, a circle is usually a tour, and it is finite.

%
\begin{prop}
\label{manydirections}
Let two finite $\Pi$-paths $\xi$ and $\eta$. If $\xi$
has many directions near a point $x\in\mathbb{R} ^d$, and its support
is included in that of $\eta$, that is, $\mathcal{L} _{\xi} \subset
\mathcal{L} _ {\eta}$, then $x \in\eta$.
\end{prop}

\begin{pf}
By Theorem~\ref{Kendall}, the finite $\Pi$-path $\eta$ has finite
Euclidean length. Moreover, $\mathcal{L} _{\xi} \subset\mathcal{L} _
{\eta}$, so that $\eta$ is a $\mathcal{L} _{\xi}$-tour in $\mathbb{R}
^d$. Since $\xi$ has many directions near~$x$, there is no finite
$\mathcal{L} _{\xi}$-tour in $\mathbb{R} ^d \setminus B(x, \varepsilon
)$, for any $\varepsilon> 0$. So that $B(x, \varepsilon) \cap\eta
\neq\varnothing$. Any finite $\Pi$-path is closed in $\mathbb{R} ^d$,
hence $x \in\eta$.
\end{pf}

%
\begin{lem}
\label{tour}
Let $ d \geq3$.

Let $l$ be a prescribed line independent of $\Pi$.
Almost surely, for all $x \in l$, for any $y \notin l$, all geodesics
$g_{xy}$ have many directions near $x$.
\end{lem}

\begin{pf}
It is sufficient to show that the claim holds for all $x \in l \cap
B(0,R)$, for all $R > 0$. It is also sufficient to work with a fixed
$\varepsilon$, and prove that for any $y \notin l^{\varepsilon}$ the
$\varepsilon$-widening of $l$, for any geodesic $g_{xy}$, all $\mathcal
{L} _{g_{xy}}$-tours on $\mathbb{R} ^d \setminus B(x, \varepsilon)$
have infinite length.

Let $r _n = r _0 / \alpha^n$ for $r_0 = \varepsilon$ and an $\alpha>
1$ big enough to be determined later. We consider one-dimensional
internal $\frac{r _n}{4} $-nets $N_n$ of $B(0, R) \cap l$.\vspace*{1pt}

In particular, there are at most $\frac{4R}{r _n } + 1$ points in $N_n$.

Now, for any point $x_n$ in $N_n$, we may build a set of nested balls
$ \{ B_i \} _{i\leq n} $ depending on $x_n$, with the following properties:
\begin{itemize}
\item{The ball $B_n$ is centred on $x_n$.}
\item{All balls $B_i = B(x_i, r _i)$ are centred on a point $x_i$ in
the net $N_i$.}
\item{The centres lie deep inside the previous ball, that is $x_i \in
B(x_{i-1}, \frac{r _{i-1}}{4})$. This stems from the fact that
$N_{i-1}$ is an $\frac{r _{i-1}}{4} $-net.}
\end{itemize}
We also use the notation of Lemma~\ref{nested_balls} for $v_i$ and
$\mathcal{V} _i$. We have to choose a good speed $v_0$.

We then show that the nested balls have the following properties, with
a $T$ to be determined later, except on an event of probability $o (
r_n^{-1} ) $:
\begin{longlist}[2.]
\item[1.] There are at least $\frac{2n}{5}$ balls $B_i$ that
no line faster than $v_{i-1}$ hits.

\item[2.] There are at least $\frac{5n}{6}$ indices $i$ such
that the time diameter of $B(x_i, \frac{7 r_i}{8}) \setminus B(x_i,
\frac{r _i}{2} ) $ on $\Pi_{[B_i]\setminus [B_{i+1}]}$ is less than $T
\alpha^{-i\vafrac{\gamma-d}{\gamma- 1}}$.

\item[3.] There are at least $\frac{5n}{6}$ balls $B_i$ that
are hit by at most $\tau$ lines of speed between $v_i$ and $v_{i-1}$,
for some fixed $\tau$. Moreover, all those lines have independent,
uniformly random directions.
\end{longlist}

The $T$ will yield the $v_0$ we need to continue the proof, that is,
%
\begin{eqnarray}
\label{v0} v_0 & =& \frac{3r_0}{8T}.
\end{eqnarray}

We deal with these properties in reverse order, and first consider
property~3.

Since the sets of lines $\mathcal{V} _i$ are disjoint, the events that
more than $\tau$ lines in $\mathcal{V} _i$ hit $B_i$ are mutually
independent (in $i$), as are their directions. Moreover, since the
$B_i$ are balls, the direction of the lines of the isotropic Poisson
line process that hit it are uniformly random.

So that the number $X_n$ of such events is a binomial random variable
$B(n, q)$ where $q$ is the probability of a single event. The Chernoff
bound for a binomial is, for $\delta> q$,
\begin{eqnarray*}
\mathbb{P} [X_n > \delta n ] & \leq&\biggl( \biggl(
\frac{q}{\delta} \biggr) ^{\delta} \biggl(\frac{1-q}{1 - \delta}
\biggr)^{1 - \delta} \biggr) ^n.
\end{eqnarray*}
Here, it suffices to take $\delta= \frac{1}6$. For $\alpha$ big
enough, if, for example, we take $q \leq\alpha^{-7}$, the above bound
is negligible with respect to $\alpha^{-n}$.

Now the probability $q$ of a single event is the probability that a
Poisson random variable with parameter not depending on $i$ is bigger
than $\tau$. We then just choose $\tau$ so that $q \leq\alpha^{-7}$.

Now for property~2.

Since the sets of lines $[B_i]\setminus[B_{i+1}]$ are disjoint, each of
the events that the time diameter of $B(x_i, \frac{7 r_i}{8}) \setminus
B(x_i, \frac{r _i}{2} ) $ on $\Pi_{[B_i]\setminus [B_{i+1}]}$ is less
than $T \alpha^{-i\vafrac{\gamma-d}{\gamma- 1}}$ are independent.
Arguing as for property~3, we merely have to choose $T$
such that the probability $q$ of a single event is less than $\alpha^{-7}$.

Now, by scaling, the sets $A_i = B(x_i, \frac{7 r_i}{8}) \setminus
B(x_i, \frac{r _i}{2} ) $ may all be covered by the same number
$\mathring{N}$ of balls with radius $\frac{r_i}{8}$, and centres in
$A_i$. Moreover, $A_i$ is connected. Now the $\frac{r_i}{8}$-widening
of $A_i$ is included in $B_i$, and since $x_{i+1} \in B(x_{i}, \frac{r
_{i}}{4})$, we know that $x_{i+1}$ is at distance at least $\frac{1}8$
of this widening. So that the maximum solid angle of $B_{i+1}$ viewed
from any point of the widening is a decreasing function of $\alpha$,
hitting zero when $ \alpha$ goes to infinity. Thus, for $\alpha$ big
enough, we may apply Lemma~\ref{Tform} with $\Omega< \Omega_{\max}$.
Replacing $\varepsilon$ in bound \eqref{td} with $1 / (\varepsilon
_{\max} \alpha^{-7})$ to ensure $q \leq\alpha^{-7}$, we find %
that property~2 is ensured if we choose
\begin{eqnarray*}
T & =& \mathring{N} T_1 \biggl(\frac{r_0}{8} \biggr)
^{\vafrac{\gamma- d}{\gamma-1} } \bigl( \ln\bigl(\alpha^7 /\varepsilon
_{\max}\bigr) \bigr) ^{\afrac{1}{\gamma-1} }
\\
& =& O \bigl( \ln^{\afrac{1}{\gamma-1} } \alpha\bigr).
\end{eqnarray*}

This choice of $T$ ensures that $v_0 = O ( \ln^{- \afrac{1}{\gamma-1} }
\alpha )$.

As for property~1, we apply Lemma~\ref{nested_balls}.

We substitute $\delta= \frac{3}{5} $ in equation \eqref{res_tower},
and use $p = \alpha^{1-d} \leq\alpha^{-2}$---since $d\geq3$---and
$\lambda = O(v_0^{1 - \gamma}) = O(\ln \alpha)$. We obtain
\begin{eqnarray*}
&& \mathbb{P} \biggl[ \mbox{For less than }\frac{2n}{5} \mbox{ indices }i,
\mbox{ no line faster than }v_{i-1}\mbox{ hits the ball
}B_i \biggr]
\\
&&\qquad \leq\bigl(c \ln(\alpha) \alpha^{-\sfrac{6}5}\bigr)^n
\\
&&\qquad  = o\bigl(\alpha^{-n}\bigr).
\end{eqnarray*}

Let us now consider a geodesic $g_{xy}$ from a point $y$ outside $B_0$
to a point $x$ inside $B_n$. Since the balls are nested, they are all
crossed by the geodesic. Even the sets $A_i = B(x_i, \frac{7 r_i}{8})
\setminus B(x_i, \frac{r _i}{2} )$ introduced in the proof of property~2 are all crossed, that is we must pass from the spherical
boundary $S(x_i, \frac{7 r_i}{8})$ to the smaller boundary $S(x_i, \frac
{r _i}{2} )$.

If the three properties are satisfied, then there are at least $\frac
{n}{15} $ indices $i$ for which the conditions are simultaneously
satisfied. For such an $i$, by property~2 the geodesic will go
from $S(x_i, \frac{7 r_i}{8})$ to $S(x_i, \frac{r _i}{2} )$ in time at
most $T \alpha^{-i\vafrac{\gamma-d}{\gamma- 1}} $. Since the two
boundaries are $ \frac{3 r_i}{8} $ apart, the geodesic must use a line
with speed more than $ \frac{3r_i}{8T \alpha^{-i(\gamma-d)/(\gamma-
1)}} = v_0 \alpha^{-i\vafrac{d-1}{\gamma-1} } = v_i $ within $A_i$.
Now,\vspace*{1pt} $A_i \subset B_i$ and, by property~1, there is no line
faster than $v_{i-1}$ in $B_i$. So that the geodesic must use a line in
$\mathcal{V} _i \cap[B_i]$, whose cardinal is less than $\tau$ by
property~3. Since the $v_i$ are disjoint, we have proved
that the geodesic must use at least $\frac{n}{15}$ distinct lines among
a set of size at most $\tau n$ of lines with uniformly random direction.

We may then apply Lemma~\ref{angles}.

We take $1 > \eta> 1 / (d-1)$ and $\delta= 1/20$. With probability $1
- O(n^{-cn})$, all geodesics from a point outside $B_0$ to a point
inside $B_n$ must use at least $n / 20$ lines (depending on the
geodesic), each with a direction in a different ball of a $n^{-\eta
}$-net of the projective space. Since there are at most a constant
number $c_1$ of points in a $n^{-\eta}$-net at distance less than $3
n^{-\eta}$ from a given point, we may choose $n / (20 c_1)$ of those
lines, so that each pair of them makes an angle at least $n^{-\eta}$.

Let us consider a fixed $r_0 = \varepsilon$. Since $d(x_i, x_{i-1})
\leq\frac{r_{i-1}}{4} $ and $x \in B(x_n, r_n) $, and since $r_i =
\varepsilon\alpha^{-i} =o(\varepsilon n^{-\eta})$ for $i > c_2 \ln
n$, we know that $B_i \subset B(x, \varepsilon n^{-\eta}/3)$ for $i >
c_2 \ln n $. Among our $n/(20 c_1)$ lines, at most $\tau$ hit each
$B_i$. Hence, up to removing $c_3 \ln n$ of our $n/(20 c_1)$ lines, all
those lines hit the ball $B(x, \varepsilon n^{-\eta}/3)$. So that, by
elementary geometry illustrated in Figure~\ref{ball}, for any two lines
$l_1$ and $l_2$ in our collection, no point of $l_1 \setminus B(x,
\varepsilon)$ is closer to a point of $l_2 \setminus B(x, \varepsilon
)$ than $\varepsilon n^{-\eta} / 4 $.

We have thus proved that for any geodesic $g_{xy}$ with $x \in B_n$ and
$y\notin B_0$, any $\mathcal{L} _{g_{xy}}$-tour in $\mathbb{R} ^d
\setminus B(x, \varepsilon)$ must contain $c_4 n$ points that are $c_5
n^{-\eta}$ apart pairwise. Hence, it has length at least $c_6 n^{1 -
\eta}$, going to infinity with $n$.

%
\begin{figure}

\includegraphics{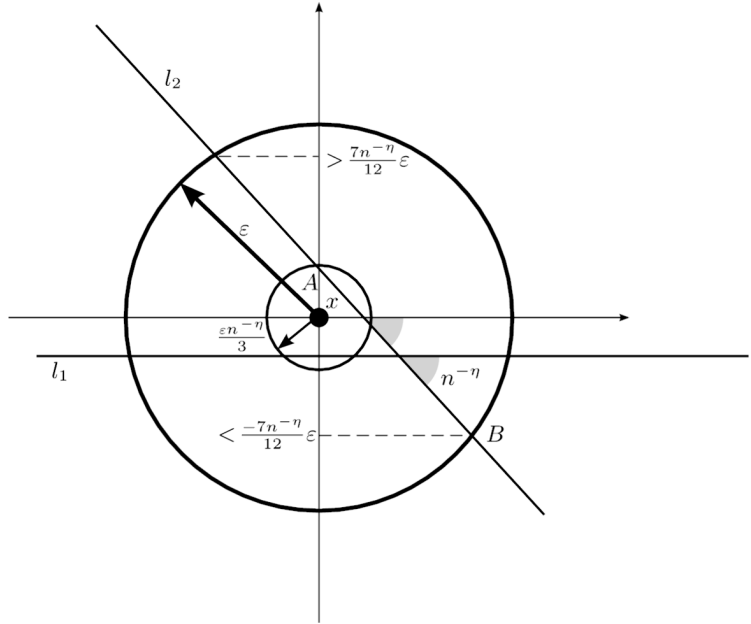}

\caption{We project on the vector space spanned by $l_1$ and $l_2$,
with $l_1$ parallel to the first coordinate and $x$ as the origin. All
points of $l_1$ then have second coordinate at most $n^{-\eta}
\varepsilon/3$ in absolute value. On the other hand, let $A$ be a
point in $l_2 \cap B(x,n^{-\eta} \varepsilon/3)$ and $B$ a point in
$l_2 \setminus B(x, \varepsilon)$. Then\vspace*{1pt} the vector $\protect
\overrightarrow{AB}$ has absolute second coordinate at least
$(\varepsilon- n^{-\eta} \varepsilon/3) \sin(n^{-\eta}) > \frac
{11}{12} \varepsilon n^{-\eta}$ for large $n$, so that $B$ has absolute
second coordinate at least $\frac{7}{12} \varepsilon n^{-\eta}$.}
\label{ball}
\end{figure}

All those properties were obtained except on a set of probability
$o(\alpha^{-n})$. Since the net $N_n$ has cardinal of order $\alpha
^{n}$, this is true simultaneously for all sets of nested balls built
on all the points $x_n$ in $N_n$, except on an event of probability
$\varepsilon_n$ going to zero.

Letting $n$ go to infinity yields the lemma.
\end{pf}

\subsection{Geodesics use the same lines}
\label{samelines}

%
\begin{teo}
\label{sameseg}
For all $x$ and $y$ in $\mathbb{R} ^d$, almost surely all geodesics
$g_{xy}$ are supported on exactly the same lines.
\end{teo}

\begin{pf}
A similar result is used during the proof of Theorem $4.4$ in
\citeauthor{Kendall2014}'s (\citeyear{Kendall2014}) article.
The idea of the proof is to ``slow down'' the lines not used by a
specific geodesic. Then all geodesics that use one of those lines
become slower and are no longer geodesics. The set of slowed speeds has
infinite measure relative to the initial speeds.

Almost surely, there exists a finite $R$ such that all geodesics from
$x$ to $y$ are included in $B(0, R)$, by Theorem $2.6$ of Kendall.
Almost surely, for any $R$, there are countably many lines intersecting
$B(0,R)$, all with different speeds $v_1 > v_2 > \cdots$ (the
information here is that there is a measurable order-preserving
bijection from $\mathbb{N} $ to the inverse line speeds). We then write
$l_i$ for the line with speed $v_i$.

Let $i \in\mathbb{N} $ and let us fix all the other speeds $v_j$ for
$j\neq i$. Let us suppose there is a speed $v^{=}_i$ between $v_{i-1}$
and $v_{i+1}$ such that if $v_i = v_i^{=}$, then there are at least two
geodesics $g_{xy}^{i}$ and $g_{xy}^{\neg i}$ from $x$ to $y$, one of
which is supported by $l_i$, and the other not. That is $l_i \sqsubset
g_{xy}^{i}$ and $l_i \not\sqsubset g_{xy}^{\neg i} $.

Now, the time length of a path not supported by $l_i$ does not depend
on $v_i$. So that $g_{xy}^{\neg i}$ is the fastest of these paths for
all $v_i$, with constant time. On the other hand, by decomposition
\eqref{ll}, the time length of $g_{xy}^{i}$ is decreasing in $v_i$.
Hence, if $v_i > v_i^{=}$, all geodesics from $x$ to $y$ are supported
by $l_i$. Conversely, if $v_i < v_i^{=}$, no path supported by $l_i$ is
a geodesic from $x$ to $y$. Indeed, such a path would be as fast as
$g_{xy}^{\neg i}$ at speed $v_i$, hence faster at speed $v_i^{=}$.

So that, when the other speeds are fixed, there is at most one value
$v_i^{=}$ of $v_i$ such that $l_i$ is in the support of some geodesic
from $x$ to $y$ and not in the support of another such geodesic.

Now, we may disintegrate the measures on the line speeds. For almost
all fixed line speeds for all lines except $l_i$, the measure $\mu_i$
for $v_i$ has a density (namely, it is proportional to $v_i^{-\gamma}
\mathbf{1}_{v_{i-1} \leq v_i \leq v_{i+1}}$). Hence, $\mu_i( \{ v_i^{=}
\}) = 0$. So that almost surely, either $l_i \sqsubset g$ for all
geodesics $g$ from $x$ to $y$, or $l_i \not\sqsubset g$ for any such
geodesic. This is true for all lines $l_i$, completing the proof.
\end{pf}

\subsection{Almost sure uniqueness}
\label{asunique}

We may now state this section's main result.

%
\begin{teo}
\label{unique_geodesics}
In dimension $d \geq2$, for all $x$ and $y$ in $\mathbb{R} ^d$,
almost surely there is a unique geodesic $g_{xy}$ from $x$ to $y$.
\end{teo}

\begin{pf}
\citet{Kendall2014} proved the dimension $2$ case in Theorem
$4.4$. We then assume from here on that $d \geq3$.

Consider any line $l \in\Pi$. We apply a Palm distribution argument:
since $\Pi$ is a Poisson process, $\Pi\setminus\{l\}$ still has the
same distribution as $\Pi$. So that almost surely, by Lemma~\ref{tour},
all geodesics on $\Pi\setminus \{ l \} $ with an endpoint $z\in l$
have many directions near $z$. Since the lines are countable, this is
true for all lines simultaneously.

The rigorous way of writing the former paragraph is through
Slivnyak--Mecke formula [originally proved by \citet{Slivnyak};
see, e.g.,\break \citeauthor{Moller}'s (\citeyear{Moller}) book, Theorem
3.2, for a modern treatment]. A Poisson point process $X$ with
intensity $\mu$ takes value in the set of locally finite point
configurations $N_{lf}$. For any nonnegative measurable function $h$
on $\mathbb{R} ^d \times N_{lf}$, we have
\begin{eqnarray*}
\mathbb{E} \biggl[ \sum_{x \in X} h\bigl(x, X \setminus
\{x\}\bigr) \biggr] & =& \int_{\mathbb{R} ^d} \mathbb{E} \bigl[ h(x, X)
\bigr] \,\mathrm{d}\mu(x).
\end{eqnarray*}
We apply the formula with the underlying point process for $\Pi$, so
that $x = (l,v)$ are marked lines. The function $h$ is the indicator function
\[
h\bigl((l,v), \Pi\bigr) = \cases{ 1, &\quad if $\exists z\in l$,
$y\notin l$,
$g_{zy}$ $\Pi$-geodesic such that
\cr
& \quad$g_{zy}$ does
not have many directions near $z$,
\cr
0, &\quad otherwise.}
\]

Lemma~\ref{tour} then ensures that the expectation in the integrand in
the right-hand side is uniformly zero, so that, almost surely, for all
$l \in\Pi$, all geodesics on $\Pi\setminus \{ l \} $ with an endpoint
$z\in l$ have many directions near $z$.

Let $g_{xy}$ be a geodesic from $x$ to $y$. Let $\{z_i\}$ be the set
of endpoints of all the segments of the geodesic $g_{xy}$ except $x$
and $y$. For any $z_i$, since it is a segment endpoint (on line $l$),
there is a $u$ such that $g_{uz_i}$ or $g_{z_i u}$ is included in
$g_{xy}$, and such that $l$ is not in the support of this sub-geodesic.
Hence, this sub-geodesic is also a geodesic in $\Pi\setminus \{ l \}
$, and has many directions near $z_i$. {A fortiori}, $g_{xy}$ has
many directions near $z_i$. Moreover, by Lemma~\ref{sameseg} almost
surely all other geodesics from $x$ to $y$ have the same support.
Hence, all these geodesics include all the segment endpoints of all the
geodesics from $x$ to $y$.

Now\vspace*{1pt} two geodesics from $x$ to $y$ must pass through the same points in
the same order:
indeed if $g^2_{xy} = g^2_{xa}g^2_{ab}g^2_{by}$ and $g^2_{xy} =
g^2_{xb}g^2_{ba}g^2_{ay}$, then $g^1_{xa}g^2_{ay}$ or
$g^2_{xb}g^1_{by}$ would be shorter than both.

So that all geodesics from $x$ to $y$ pass through their segment
endpoints in the same order, so they are the same.
\end{pf}

A byproduct of the proof is the following remark.
%
\begin{cor}
\label{strengthening}
Let $d \geq3$. Almost surely, for any point $x$ not on a line of~$\Pi
$, that is $x \notin\mathcal{S} $, all geodesics containing $x$ have
many directions near $x$.
\end{cor}

\begin{pf}
We use the step in the former proof, that almost surely, for all $l
\in\Pi$, all geodesics on $\Pi\setminus \{ l \} $ with an endpoint
$z\in l$ have many directions near $z$.

For any $\varepsilon$, since $x \notin\mathcal{S} $, the geodesic to
$x$ will leave a line of $\Pi$ at a point $z \in B(x, \varepsilon/2)$.
Now the restriction of the geodesic to $g_{zx}$ has many directions
near~$z$. Since $B(z, \varepsilon/ 2 ) \subset B(x, \varepsilon)$,
$g$ has also many directions near $x$.
\end{pf}

\section{Geodesic length has finite expectation}
\label{expect_geo}

This is property~3 of a SIRSN.

%
\begin{teo}
\label{teo_len}
Let $x,y \in\mathbb{R} ^d$. Then the Euclidean length of the $\Pi
$-geodesic between $x$ and $y$ has finite expectation
\begin{eqnarray*}
\mathbb{E} [L_{xy} ] & <& \infty.
\end{eqnarray*}
\end{teo}

\begin{pf}
For any $r$, any $\Pi$-path containing $x$ whose Euclidean length is
more than $r$ has to intersect $ B(x,r)$ on a length at least $r$.
Hence, it must spend time at least $r / v_{\max}(B(x,r))$ in that ball.
If the geodesic $g_{xy}$ has Euclidean length $L_{xy}$ and time length
$T_{xy}$, we thus obtain the constraint
%
\begin{eqnarray}
\label{constraint} v_{\max}\bigl(B(x,r)\bigr) & \geq& r / T_{xy}
\qquad\mbox{for all }r \leq L_{xy}.
\end{eqnarray}

By equation \eqref{td} of Theorem~\ref{Tform}, there is a $T$ such that
with probability at least $1 - 2^{-(n+1)}$
%
\begin{eqnarray}
\label{n_time} T_{xy} & \leq&(n+1)^{\afrac{1}{\gamma- 1} } T \,\hat{=}\, T_n.
\end{eqnarray}


We now consider a collection of radii
\begin{eqnarray*}
r_l & =& 2^l r_0
\end{eqnarray*}
for $l$ between $0$ and $m$, with $r_0$ and $m$ to be chosen later,
possibly depending on $n$, but independent of $\Pi$.

If $L_{xy} > r_m$, then $L_{xy} > r_l$ for all $l \in[0,m]$, so that
constraint \eqref{constraint} must be satisfied for each $r_l$. In
particular, on the event \eqref{n_time}, the following constraint is
satisfied for all $l \in[0,m]$:
%
\begin{equation}
\label{constraint2} v_{\max}\bigl(B(x,r_l)\bigr) \geq
r_l / T_n \,\hat{=}\, v_l.
\end{equation}
We again drop the dependence on $n$ in the notation of $v_l$, and $v_l$
is independent of~$\Pi$.

These constraints are simultaneously satisfied if and only if there is
a strictly increasing sequence of $(1+k)$ integers $0=l_0, l_1, \dots,
l_k$ between $0$ and $m$ such that
\begin{eqnarray*}
v_{l_{i+1}} &>& v_{\max}\bigl(B(x,r_{l_i})\bigr)  \geq
v_{l_{i+1} - 1}\qquad\mbox{if }i < k,
\\
v_{\max}\bigl(B(x,r_{l_k})\bigr) & \geq& v_m.
\end{eqnarray*}

We define the following events for all $i\in[0,m]$, with the
conventions $l_{k+1} = m+1$ and $r_{l_{-1}} = 0$:
\begin{eqnarray*}
A_i & =& \bigl\{ v_{\max}\bigl(B(x,r_{l_i})\bigr)
\geq v_{l_{i+1} - 1} \bigr\},
\\
B_i & =& \bigl\{ v_{l_{i+1}} > v_{\max}
\bigl(B(x,r_{l_i})\bigr) \geq v_{l_{i+1} - 1} \bigr\},
\\
D_i & =& \bigl\{ v_{\max} \bigl(\bigl[B(x,r_{l_i})
\bigr]\setminus\bigl[B(x,r_{l_{i-1}})\bigr] \bigr) \geq v_{l_{i+1} - 1}
\bigr\}.
\end{eqnarray*}

Notice that the events do not depend only on $i$, but on the whole
sequence of $l_i$. Notation is easier this way. In particular, the
former paragraph reads
%
\begin{eqnarray}
\label{borneLxy} \qquad\mathbb{P} [ L_{xy} > r_m |
T_{xy} \leq T_n ] & \leq&\sum_{k=0}^m
\sum_{0=l_0 < \cdots< l_k \leq m} \mathbb{P} [B_0 \cap
B_1 \cap\cdots\cap B_{k-1} \cap A_k ].
\end{eqnarray}

Let us consider the filtration generated by the lines intersecting
$B(x,r)$ for increasing $r$, that is $\mathcal{F} _r = \sigma ( \Pi
_{[B(x,r)]} ) $. Then $A_i$, $B_i$ and $D_i$ are all $\mathcal{F}
_{r_{l_i}}$-measurable, and $D_i$ is independent of $\mathcal{F}
_{r_{l_{i-1}}}$.

Moreover, the difference between $A_i$ and $D_i$ is only on the event
that there is a line faster than $v_{l_{i+1} - 1}$ hitting
$B(x,r_{l_{i-1}}) $. Since $l_{i+1} - 1 \geq l_i$, this never happens
under $B_{i-1}$. So that $A_i \cap B_{i-1} = D_i \cap B_{i-1}$. From
this, we deduce
\begin{eqnarray}\label{Pai}
\mathbb{P} [B_i | B_0, \dots, B_{i-1} ] & \leq&\mathbb{P} [A_i | B_0, \dots, B_{i-1} ]
\qquad\mbox{since }B_i \subset A_i\nonumber
\\[-2pt]
& =& \mathbb{P} [D_i | B_0, \dots, B_{i-1} ]
\qquad\mbox{since }A_i \cap B_{i-1} = D_i
\cap B_{i-1}
\nonumber\\[-9pt]\\[-9pt]\nonumber
& =& \mathbb{P} [D_i ]\qquad\mbox{since }D_i\mbox{
indep. of }B_j\mbox{ for }j < i
\\[-2pt]
& \leq&\mathbb{P} [A_i ]\qquad\mbox{since }D_i \subset
A_i. \nonumber
\end{eqnarray}

Recall that the number of lines faster than $v$ hitting a ball of
radius $r$ is a Poisson variable with parameter $c r^{d-1} v^{-(\gamma
-1)}$, with $c = \mu_d([B(0,1)])$. We may then compute
\begin{eqnarray}\label{unP}
\mathbb{P} [A_i ] & =& 1 - \exp\bigl( - c r_{l_i}^{d-1}
v_{l_{i+1} - 1} ^{-(\gamma - 1)} \bigr)\nonumber
\\[-2pt]
& \leq& c r_{l_i}^{d-1} v_{l_{i+1} - 1} ^{-(\gamma - 1)}\nonumber
\\[-2pt]
& =& c r_{l_i}^{d-1} 2^{-(\gamma- 1)(l_{i+1} - l_i -1)} v_{l_i}
^{-(\gamma - 1)}\nonumber
\\[-2pt]
& =& 2^{-(\gamma- 1)(l_{i+1} - l_i -1)} c r_{l_i}^{d- \gamma}
T_n^{\gamma- 1}
\\[-2pt]
& \leq&2^{-(\gamma-1)(l_{i+1} - l_i -1)} c r_{0}^{d- \gamma}
T_n^{\gamma- 1}\qquad\mbox{since }r_0 \leq
r_{l_i}\nonumber
\\[-2pt]
& =& 2^{-(\gamma-1)(l_{i+1} - l_i)} p(r_0) \qquad\mbox{with}\nonumber
\\[-2pt]
p(r_0)  &=& 2^{\gamma-1} c r_{0}^{d- \gamma}
T_n^{\gamma- 1}.\nonumber
\end{eqnarray}

Recalling the convention $l_{k+1} = m+1$, we obtain
\begin{eqnarray*}
&& \mathbb{P} [B_0 \cap B_1 \cap\cdots
\cap B_{k-1} \cap A_k ]
\\[-2pt]
&&\qquad  \leq\mathbb{P} [B_0 ] \mathbb{P} [B_1 |
B_0 ] \cdots\mathbb{P} [A_k | B_0,
B_1,\dots, B_{k-1} ]
\\[-2pt]
&&\qquad  \leq\prod_{i=0}^{k} \mathbb{P}
[A_i ]\qquad\mbox{by \eqref{Pai}}
\\[-2pt]
&&\qquad  \leq\prod_{i=0}^{k} 2^{-(\gamma-1)(l_{i+1} - l_i)}
p(r_0)\qquad\mbox{by \eqref{unP}}
\\[-2pt]
&&\qquad  = 2^{-(\gamma- 1)(m+1)} p(r_0)^{k+1}.
\end{eqnarray*}

Substituting into bound \eqref{borneLxy}, we obtain
\begin{eqnarray}
\mathbb{P} [ L_{xy} > r_m | T_{xy} \leq
T_n ] & \leq&\sum_{k=0}^m
\sum_{0=l_0 < \cdots< l_k \leq m} 2^{-(\gamma- 1)(m+1)} p(r_0)^{k+1}\nonumber
\\
\label{sum} & =& \sum_{k=0}^m \pmatrix{m
\cr k} 2^{- (m+1)(\gamma-1)} p(r_0)^{k+1}
\\
& =& 2^{- (m+1)(\gamma-1)} p(r_0) \bigl(1 + p(r_0)
\bigr)^m\nonumber
\\
& \leq&\bigl( 2^{-(\gamma-1)} \bigl(1 + p(r_0)\bigr) \bigr)
^{m+1}. \label{reborne}
\end{eqnarray}

Let $\kappa< \gamma-1$. We now choose our free parameters:
\begin{eqnarray*}
r_0 & =& \biggl( \frac{2^{\gamma- 1} c T^{\gamma- 1}}{2^{\gamma-1 -
\kappa} - 1} \biggr) ^{\afrac{1}{\gamma- d}},
\end{eqnarray*}
so that
\begin{eqnarray*}
1 + p(r_0) & =& 2^{\gamma-1 - \kappa},
\\
m & =& \bigl\lfloor(n+1) / \kappa\bigr\rfloor.
\end{eqnarray*}

Substituting into bound \eqref{reborne}, we get
\begin{eqnarray*}
\mathbb{P} [ L_{xy} > r_m | T_{xy} \leq
T_n ] & \leq&2^{-(m+1) \kappa}
\\
& \leq&2^{-(n+1)}.
\end{eqnarray*}

Since $\mathbb{P} [ T_{xy} > T_n ] \leq2^{-(n+1)}$, we have thus
proved that
with probability at least $1 - 2^{-n}$, the Euclidean length is bounded by
\begin{eqnarray*}
L_{xy} & \leq& r_m
\\
& =& 2^{(n+1) / \kappa} (n+1)^{\afrac{1}{\gamma- d} } C,
\end{eqnarray*}
where $C$ is a positive constant depending on $T$ and $\kappa$, but
not on $n$.

Thus, $L_{xy}$ has a $\delta$-moment for all $\delta< \kappa$. Indeed,
\begin{eqnarray*}
\mathbb{E} \bigl[ L_{xy}^{\delta} \bigr] & \leq&\sum
_{n=0}^{\infty} 2^{-n} \bigl( 2^{(n+1) / \kappa}
(n+1)^{\afrac{1}{\gamma- d} } C \bigr)^{\delta}
\\
& \leq&\sum_{n=0}^{\infty} 2^{(\delta/ \kappa- 1) n } O
\bigl( n^{\afrac{\delta}{\gamma-d} } \bigr)
\\
& <& \infty \qquad\mbox{if }\delta< \kappa.
\end{eqnarray*}
Since $\kappa$ is only constrained by $\kappa < \gamma- 1$, we have
a $\delta$-moment for all $\delta< \gamma- 1$. Since $\gamma> d
\geq2$, this completes the proof.
\end{pf}

%
\begin{rem}
The main part of the proof is really just saying that the
Euclidean diameter of a $\Pi$-ball has finite mean. Together with
Theorem~\ref{Tform}, this implies that the random metric space
generated by the $\Pi$-length is almost surely homeomorphic to $\mathbb
{R} ^d$.

Notice that for $\kappa$ close to $\gamma- 1$, the main term in
the sum \eqref{sum} is the one for $k=0$. In other words, it is easier
to have one single extremely fast line close to $x$ than to have many
successively faster lines, if we want to be abnormally fast at each
distance of $x$.

A first way to improve on the moments of $L_{xy}$ starts with
noticing that if we use a line to go far away very fast, we need to use
another line to come back, since a geodesic never crosses itself.

More precisely, the author conjectures that the structure
yielding long geodesics with highest probability is the following: a
line with speed at least $2r/T_{xy}$ hits $B(x, \varepsilon)$, another
hits $B(y, \varepsilon)$ and they hit a common $\varepsilon$-ball at
distance~$r$. Since the two first events have probability in
$r^{-(\gamma- 1)}$ and the other in $r^{-(d - 1)}$, we would conclude
that $\mathbb{E} [ L_{xy}^{2 \gamma+ d - 3} ] = \infty$, but $\mathbb
{E} [ L_{xy}^\delta ] < \infty$ for all $\delta< 2 \gamma+ d - 3$.
\end{rem}

\section{Finite intensity of long-distance network}
\label{ldn}

We now turn to property~4 of a SIRSN, that is their key property.

Intuitively, this means that the SIRSN contains ``highways.'' If we
look at all the geodesics simultaneously, truncating each geodesic by
deleting balls around its endpoints, their total length in each compact
set is finite: the geodesics largely re-use the same segments in each
region. Contrast with the Euclidean case where the whole space is used.

In our context, we have to prove that the intensity $p(1)$ of the
following long-distance network $\mathcal{F} $ is finite: let $ \{ \Xi
_n, n\in\mathbb{N} ^* \} $ be a collection of Poisson processes on
$\mathbb{R} ^d$ with intensity $n$ times Lebesgue, all independent from
$\Pi$, and coupled so that $\Xi_n \subset\Xi_{n+1}$. Write $\Xi=
\bigcup_{n\in\mathbb{N} ^*} \Xi_n$. Then
\begin{eqnarray*}
\mathcal{F} & =& \bigcup_{x,y \in\Xi} \bigl( g_{xy}
\setminus\bigl(B(x, 1) \cup B(y,1) \bigr)\bigr).
\end{eqnarray*}

Notice that almost surely $g_{xy}$ is unique for all $x$ and $y$, since
the dense point set $\Xi$ is countable.

By translational invariance, it is enough to prove that the
intersection of $\mathcal{F} $ with a given ball has finite mean
Hausdorff measure of dimension $1$. Indeed, if $\ell= m_1(\mathcal{F}
\cap B(x,r) )$, then $\mathbb{E} [\ell] = \frac{\omega
_{d-1}r^{d-1}}{2} p(1)$. Notice that scale equivariance yields similar
results if the long-distance network $\mathcal{F} $ is defined by
removing balls of any fixed radius instead of radius $1$.

The main argument relies on the pigeon-hole principle: a geodesic
getting close to a prescribed point must use one of a very few fast
lines close to that point, and must use them again to draw away. And by
uniqueness of geodesics, two geodesics with two common points must
agree between those points.

%
\begin{teo}
\label{finite_intensity}
Let $\gamma> d \geq2$.
With the above notation, let $\ell= m_1(\mathcal{F} \cap B(0, \frac
{1}{3}) )$ be the length of the long-distance network in $B(0, \frac
{1}{3} )$. For $\varepsilon< \varepsilon_{\max}$, with probability at
least $1 - \varepsilon$, this length is less than $C (\ln(C_1 /
\varepsilon))^2$, for constants $\varepsilon_{\max}$, $C$ and $C_1$
depending only on $\gamma$ and $d$. Consequently, there is a finite
moment of exponential form: for any $\delta< \sqrt{C}$,
%
\begin{eqnarray}
\label{moment_l} \mathbb{E} \bigl[\exp(\delta\sqrt{\ell} ) \bigr] & <&
\infty.
\end{eqnarray}
In particular, $\ell$ has finite mean.
\end{teo}

\begin{pf}
Since $\Xi$ is countable, by Theorem~\ref{unique_geodesics}, almost
surely all geodesics between its points are unique. In the proof, we
only consider such geodesics and sub-geodesics, so that we assume uniqueness.

Consider the balls $B(0, \frac{1}{3} )$ and $B(0, \frac{2}{3} )$. We
call their set difference $B(0, \frac{2}{3} ) \setminus B(0, \frac
{1}{3} )$ the \emph{border}.

Now if a point $x \in B(0, \frac{2}{3})$, then $B(0, \frac{1}{3} )
\subset B(x,1)$. Hence, geodesics with an endpoint in $B(0, \frac{2}{3}
)$ make no contribution to $\ell$. We have
\begin{eqnarray*}
\ell& \leq& m_1 \biggl( \bigcup_{x,y \in\Xi\setminus B(0, \sfrac
{2}{3}) }
\biggl( g_{xy} \cap B\biggl(0, \frac{1}{3} \biggr)\biggr)
\biggr).
\end{eqnarray*}


%
\begin{figure}[b]

\includegraphics{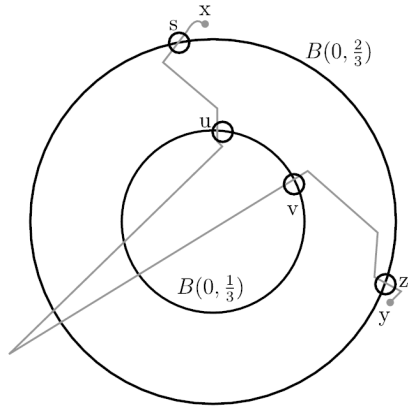}

\caption{The geodesic from $x$ to $y$ hits $B(0, \frac{2}{3} )$ for
the first time at $s$ and the last time at $z$. It hits $B(0, \frac
{1}{3} )$ for the first time at $u$ and the last time at $v$.}
\label{local_geo}
\end{figure}

Hence, geodesics $g_{xy}$ making contributions to $\ell$ are
structured in the following way, illustrated in Figure~\ref{local_geo}.
\begin{itemize}
\item{They hit $B(0, \frac{2}{3} ) $ for the first time at point $s$ on
the corresponding sphere.}
\item{Then they hit $B(0, \frac{1}{3} )$ for the first time at point
$u$ on the corresponding sphere.}
\item{Then they hit $B(0, \frac{1}{3}) $ for the last time at point $v$
on the corresponding sphere.}
\item{Then they hit $B(0, \frac{2}{3} ) $ for the last time at point
$z$ on the corresponding sphere.}
\end{itemize}
In particular, the contribution to $\ell$ is included in the
sub-geodesic $g_{uv}$.

Now, uniformly on $x$ and $y$, the time $T_{sz}$ between $s$ and $z$ is
bounded by the time diameter of $B(0, \frac{2}{3} )$. So that, by
Theorem~\ref{Tform}, with probability $1 - \varepsilon/ 2$, we have
%
\begin{eqnarray}
\label{Teps} T_{sz} & \leq& T_{\varepsilon}
\nonumber\\[-9pt]\\[-9pt]\nonumber
& \,\hat{=}&  T_{2/3, 1} \biggl( \ln\frac{2 }{\varepsilon\cdot
\varepsilon_{\max} } \biggr)^{\afrac{1}{\gamma- 1}}.
\end{eqnarray}

We\vspace*{2pt} call \emph{fast lines} the lines faster than $v_{\varepsilon} = 1 /
(6 T_{\varepsilon})$. Conversely, lines slower than $v_{\varepsilon}$
are \emph{slow lines}. We write $\mathcal{V} _{\varepsilon} = \mathcal
{S} _{v_{\varepsilon}} \cap B(0, \frac{2}{3} )$ for the intersection
of all those fast lines with the ball $B(0, \frac{2}{3} )$.

The number of fast lines hitting $B(0, \frac{2}{3} )$ is a Poisson
variable with parameter $\lambda_{\varepsilon} = \frac{\omega_{d-1}
(2/3)^{d-1}}{2 v_{\varepsilon}^{\gamma- 1}} = C_2 \ln(C_1 /
\varepsilon)$ for\vspace*{1pt} explicit constants $C_1$ and $C_2$, with $C_1 =
2/\varepsilon_{\max} \geq2$. We recall that the moment generating
function of such a Poisson variable $X$ is $\mathbb{E} [e^{t X} ] = \exp
(\lambda_{\varepsilon} (e^t - 1))$, and use Chernoff bound to get
\begin{eqnarray*}
\mathbb{P} [ X \geq C_3 \lambda_{\varepsilon} ] & \leq&
e^{\lambda_{\varepsilon} (e^t - 1)} e^{- t C_3 \lambda_{\varepsilon
}}
\\[-2pt]
& =& \exp\bigl( \lambda_{\varepsilon} (C_3 - 1 -
C_3 \ln C_3) \bigr) \qquad\mbox{with }
e^t = C_3
\\[-2pt]
& =& \biggl(\frac{C_1 }{\varepsilon} \biggr)^{C_2 (C_3 - 1 - C_3 \ln
C_3)}
\\[-2pt]
& \leq&\frac{\varepsilon}{2},
\end{eqnarray*}
with $C_3$ chosen big enough to have $C_2 (C_3 - 1 - C_3 \ln C_3) \leq
-1 $.

%
%
\begin{figure}[b]

\includegraphics{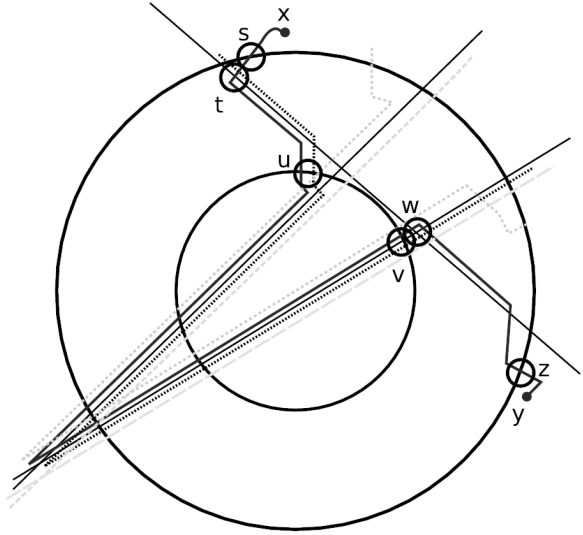}

\caption{Close parallel curves in the figure agree. They are separated
to make the figure more readable. The three thin black lines are the
fast lines. The light dashed lines and black dotted curves are a family
of geodesics $\{g^i\}$. The solid curve connecting $x$ and $y$ is the
geodesic $g_{xy}$.
The geodesic $g_{xy}$ has a common point $t$ with the black dotted
geodesic on a fast (black) line when first crossing the border, and
another $w$ when last crossing the border. Hence, they agree between
$u$ and $v$.
Any other geodesic contributing to $\ell$ would meet one of the
geodesics in the family $\{g^i\} $ in the same way.}
\label{pigeonhole}
\end{figure}

Hence, with probability at least $1 - \frac{\varepsilon}{2}$, there
are at most $C_3 \lambda_{\varepsilon}$ fast lines hitting $B(0, \frac
{2}3)$. With probability at least $1 - \varepsilon$ both this event
and $T_{sz} \leq T_{\varepsilon}$ are true. We assume both from now on.

Since $T_{sz} \leq T_{\varepsilon}$, the intersection of $g_{sz}$
with slow lines has length at most $v_{\varepsilon} / T_{\varepsilon}
= \frac{1}{6} $.

So that, since $s$ and $u$ (resp., $v$ and $z$) are at least
$\frac{1}3$ apart, the geodesic must have fast segments for length at
least $\frac{1}{6} $ between $s$ and $u$ (resp., $v$ and $z$),
that is,
\begin{eqnarray*}
m_1( \mathcal{V} _{\varepsilon} \cap g_{su}) & \geq&
\tfrac{1}{6},
\\[-2pt]
m_1( \mathcal{V} _{\varepsilon} \cap g_{vz}) & \geq&
\tfrac{1}{6}.
\end{eqnarray*}

Now
\begin{eqnarray*}
m_1(\mathcal{V} _{\varepsilon}) & \leq& \# \{ \mbox{fast lines} \}
\cdot\Diam\biggl(B\biggl(0, \frac{2}{3} \biggr) \biggr)
\\[-2pt]
& \leq&\frac{4 C_3 \lambda_{\varepsilon}}{3}
\\[-2pt]
& =& C_4 \ln\biggl( \frac{C_1}{\varepsilon} \biggr),
\end{eqnarray*}
so that\vspace*{-2pt}
\begin{eqnarray*}
\frac{ m_1\otimes m_1(\mathcal{V} _{\varepsilon} \otimes\mathcal{V}
_{\varepsilon})}{m_1 \otimes m_1 (g_{su} \cap\mathcal{V}
_{\varepsilon} \otimes g_{vz} \cap\mathcal{V} _{\varepsilon} ) } & \leq& C_5 \biggl( \ln\frac{C_1}{\varepsilon}
\biggr)^2.
\end{eqnarray*}

Hence, by the pigeon-hole property as illustrated in Figure~\ref
{pigeonhole}, we may find a maximal family $ \{ g^i \}$ of at most $
(C_5 ( \ln\frac{C_1}{\varepsilon} )^2 + 1)$ geodesics such that:
\begin{itemize}
\item{$g^i$ is a geodesic between $x^i$ and $y^i$ in $\Xi\setminus
B(0, \frac{2}{3} )$, passing through $s^i$, $u^i$, $v^i$ and $z^i$
defined as for $s,u, v$ and $z$ above.}
\item{Any geodesic $g_{xy}$ contributing to $\ell$ crosses one of the
$g^i$ when first and last crossing the border, that is, there are
points $t$ and $w$ such that $t \in g_{s^iu^i} \cap g_{su} $ and $w \in
g_{v^iz^i} \cap g_{vz} $.}
\end{itemize}
By uniqueness of geodesics, $g^i$ and $g_{xy}$ coincide on $g_{tw}$. In
particular, they coincide on $g_{uv}$. Hence, the intersection $\mathcal
{F} \cap B(0, \frac{1}{3} )$ is included in the finite number of
geodesics $g_{u^iv^i}$. We may then conclude by separating
contributions from fast and slow lines:
\begin{eqnarray*}
\ell& \leq&\sum_i m_1 (
g_{u^iv^i} )
\\
& \leq& m_1(\mathcal{V} _{\varepsilon}) + \sum
_i m_1 \bigl( g_{u^iv^i} \cap(\Pi
\setminus\mathcal{V}_{\varepsilon}) \bigr)
\\
& \leq& C_4 \ln\biggl( \frac{C_1}{\varepsilon} \biggr) +
\biggl(C_5 \biggl( \ln\frac{C_1}{\varepsilon} \biggr)^2 + 1
\biggr) \cdot\frac{1}{6}
\\
& \leq& C \biggl( \ln\frac{C_1}{\varepsilon} \biggr)^2.
\end{eqnarray*}\upqed
\end{pf}

%
\begin{teo}
\label{isSIRSN}
The network $\mathcal{N} $ made of the time geodesics is a SIRSN.
\end{teo}

\begin{pf}
Property~1 of a SIRSN is a consequence of almost sure
uniqueness of geodesics between two points, that is Theorem~\ref
{unique} in dimension at least $3$, and Theorem 4.4. in \citeauthor
{Kendall2014}'s (\citeyear{Kendall2014}) article in dimension $2$.

Property~2 of a SIRSN is because the underlying Poisson
line process is invariant by translation and rotation. As for change of
scale, the underlying Poisson line process is invariant by a
transformation where scale is multiplied by $\alpha$ and speed by
$\alpha^{(d-1)/(\gamma- 1)}$. Hence, all paths have their time length
multiplied by the same $\alpha^{(\gamma-d) / (\gamma-1)}$, so that
the geodesics are the same and $\mathcal{N} $ is invariant.

Property~3 of a SIRSN is Theorem~\ref{teo_len}.

Property~4 of a SIRSN is Theorem~\ref{finite_intensity}.
\end{pf}

\section{Conclusion}
\label{conclusion}

We have established that the improper Poisson line process with
adequate speed limits yield a SIRSN.

Along the way, a few questions have been raised. Is there an easier,
more natural way to prove uniqueness of geodesics? What are the
tightest moments of the Euclidean length of a geodesic? When can we
generalize this construction using geodesics from a random geodesic
metric space?

On a more general note, we may wonder which property of our network
$\mathcal{N} $ translate to general SIRSNs, or to SIRSNs made of
geodesics of a metric space. For example, it should be easy to show
that not all SIRSNs have the equivalent of Lemma~\ref{tour} or
Corollary~\ref{strengthening}: it is certainly not true for the binary
hierarchy model by Aldous.

We might also raise a few typical questions in stochastic geometry. Is
there only one geodesic connecting a prescribed point $x$ to infinity,
like in dimension two? In many models, infinite geodesics have an
asymptotic direction. For SIRSNs, this property looks unlikely, and
characterizing the random walk of the angle as a function of the $\Pi
$-distance to $x$ looks worthwhile. What is the law of a typical cell
in the tessellation generated by the network connecting the points of
$\Xi_1$, an intensity $1$ Poisson point process?

Finally, another somewhat tangential direction of research would be to
study more closely the properties of the random metric space. For
example, being a SIRSN entails coalescence of geodesics, a very
hyperbolic-like property. We may also draw comparisons with a
well-known random metric space such as the Brownian map
[see, e.g., \citet{LeGall2014}].

The Brownian map is a random metric space homeomorphic to the sphere
$\mathbb{S} ^2$. It has Hausdorff dimension $4$. All its geodesics
minus their endpoints is of Hausdorff dimension $1$. The cut-locus of
its geodesics starting from a given point has Hausdorff dimension $2$,
and the topology of an open continuous tree.

On the other hand, we have shown that our metric space is homeomorphic
to~$\mathbb{R} ^d$. It should be easy, by scaling arguments, to show
that its Hausdorff dimension is $(d \gamma-d) / (\gamma- d)$, which
is bigger than $d$. In particular, with $d=2$ and $\gamma= 3$, we have
the same dimensions as the Brownian map. If any geodesic can be
appropriately approximated by geodesics between points of Poisson point
processes, it should also be easy to show that the geodesics minus
their endpoints is of dimension $1$. However, the cut-locus might have
a very different behaviour.


\section*{Acknowledgements}
I would like to thank Wilfrid Kendall for introducing us to the
problem during a talk, and helpful discussions since then. Remarks from
a referee have greatly improved the presentation of the paper.


%

\printaddresses
\end{document}